\documentclass[a4paper,12pt]{article}
\usepackage[utf8]{inputenc}

\title{Two novel families of multiscale staggered patch schemes efficiently simulate large-scale, weakly damped, linear waves}

\author{
J. Divahar
\thanks{School of Mathematical Sciences,
University of Adelaide, South Australia.}
\thanks{\protect\url{https://orcid.org/0000-0002-9506-8846},
\protect\url{mailto:jdivahar@outlook.com}}
\and
A.~J. Roberts
\footnotemark[1]
\thanks{\protect\url{http://orcid.org/0000-0001-8930-1552}}
\and
Trent W. Mattner
\footnotemark[1]
\thanks{\protect\url{https://orcid.org/0000-0002-5313-5887}}
\and
J.~E. Bunder
\footnotemark[1]
\thanks{\protect\url{http://orcid.org/0000-0001-5355-2288}}
\and
Ioannis~G. Kevrekidis
\thanks{Departments of Chemical and Biomolecular Engineering \& Applied Mathematics and Statistics, Johns Hopkins University, Baltimore, Maryland, USA.
\protect\url{https://orcid.org/0000-0003-2220-3522}}
}

\usepackage{amsfonts,amsmath,amssymb,amsthm}

\IfFileExists{ajrX.sty}{% AJR style preference
  \usepackage[a5paper,margin=6.7mm,tmargin=13mm]{geometry}
  \usepackage{trackRefs}
  \AtBeginDocument{
    \setNodeShape{tri}
    \setuptodonotes{inline,color=orange!25}
    \addbibresource{ajr,bib}
  }
  \usepackage{wrapfig}

}{ % alternate style preference
  \usepackage[small, euler-digits]{eulervm} % for math
  \usepackage{palatino}  % for text
  \let\vec\mathbf
  \AtBeginDocument{
    \setNodeShape{circ}  % "circ" or "tri" to set the shapes that illustrate h,u,v nodes in the figures and the text
    \setuptodonotes{inline,color=orange!25}
  }
  
}

\renewcommand{\vec}[1]{\text{\boldmath$#1$}}
\allowdisplaybreaks

\usepackage{microtype} % https://en.wikipedia.org/wiki/Microtypography
\usepackage{needspace}
\usepackage[defaultlines=2,all]{nowidow}% seems useful to have

\usepackage[svgnames]{xcolor}
\usepackage{graphicx,subcaption}
\graphicspath{{figs/inkscape/}}
\usepackage[leftcaption]{sidecap}
\sidecaptionvpos{figure}{t} % To have the caption and the figure top-aligned

\usepackage{multirow}
\usepackage{algorithm}
\usepackage{algpseudocode}
\algrenewcommand{\algorithmiccomment}[1]{{\color{Grey}\small #1}}

\usepackage{hyperref}
\usepackage{hypcap}
\hypersetup{
  pdftex,
  backref=true,
  hyperindex=true,
  colorlinks=true,
  allcolors=teal,
  bookmarks=true,
  breaklinks=true
}

\usepackage[capitalise,nameinlink,noabbrev]{cleveref}
\crefname{equation}{}{} % to get ``(1)'' instead of ``Equation (1)''
\crefname{enumi}{}{}
\crefname{enumii}{}{}
\crefname{enumiii}{}{}
\crefname{enumiv}{}{}

\usepackage{citeSettings}
\makeatletter
\let\L@TeX@startsection\@startsection
\renewcommand{\@startsection}[6]{\L@TeX@startsection{#1}{#2}{#3}{#4}{#5}
{\raggedright #6}}
\makeatother

\usepackage{defs}
\def\xv{\vec x}

\usepackage[textwidth=3.3cm,textsize=small]{todonotes}
\usepackage{marginnote}

\begin{document}

\maketitle

\begin{abstract}
  Many multiscale wave systems exhibit macroscale emergent behaviour, for example, the fluid dynamics of floods and tsunamis.
  Resolving a large range of spatial scales typically requires a prohibitively high computational cost.
  The small dissipation in wave systems poses a significant challenge to further developing multiscale modelling methods in multiple dimensions.
  This article develops and evaluates two families of equation-free multiscale methods on novel 2D staggered patch schemes, and demonstrates the power and utility of these multiscale schemes for weakly damped linear waves.
  A detailed study of sensitivity to numerical roundoff errors establishes the robustness of developed staggered patch schemes.
  Comprehensive eigenvalue analysis over a wide range of parameters establishes the stability, accuracy, and consistency of the multiscale schemes.
  Analysis of the computational complexity shows that the measured compute times of the multiscale schemes may be \(10^5\) times smaller than the compute time for the corresponding full-domain computation.
  This work provides the essential foundation for efficient large-scale simulation of challenging nonlinear multiscale waves.
\end{abstract}

\tableofcontents

\section{Introduction}
% ++++++++++++++++++++++++++++++++++++++++++++++++++++++++++++++++++++++++++++++

\begin{figure}
  \caption{\label{fig:tSim2_gLnrWve_h_r0p1}%
  Wave height~\(h\) of time simulation from a hump with initial velocity along~\(x\), using a patch scheme (square-p4) over a staggered patch grid with \(18 \times 18\) macro-grid intervals (\(N=18\)) and each patch containing \(6 \times 6\) sub-patch micro-grid intervals (\(n=6\)), and patch ratio~\(r = 0.1\).
  %
  % Here the spatial patches are enlarged to be seeable.  % [JD] Without enlarging, it was 'seeable', only not clear.
  Here the spatial patches are enlarged for visual clarity.
  }
  \centering
  \includegraphics{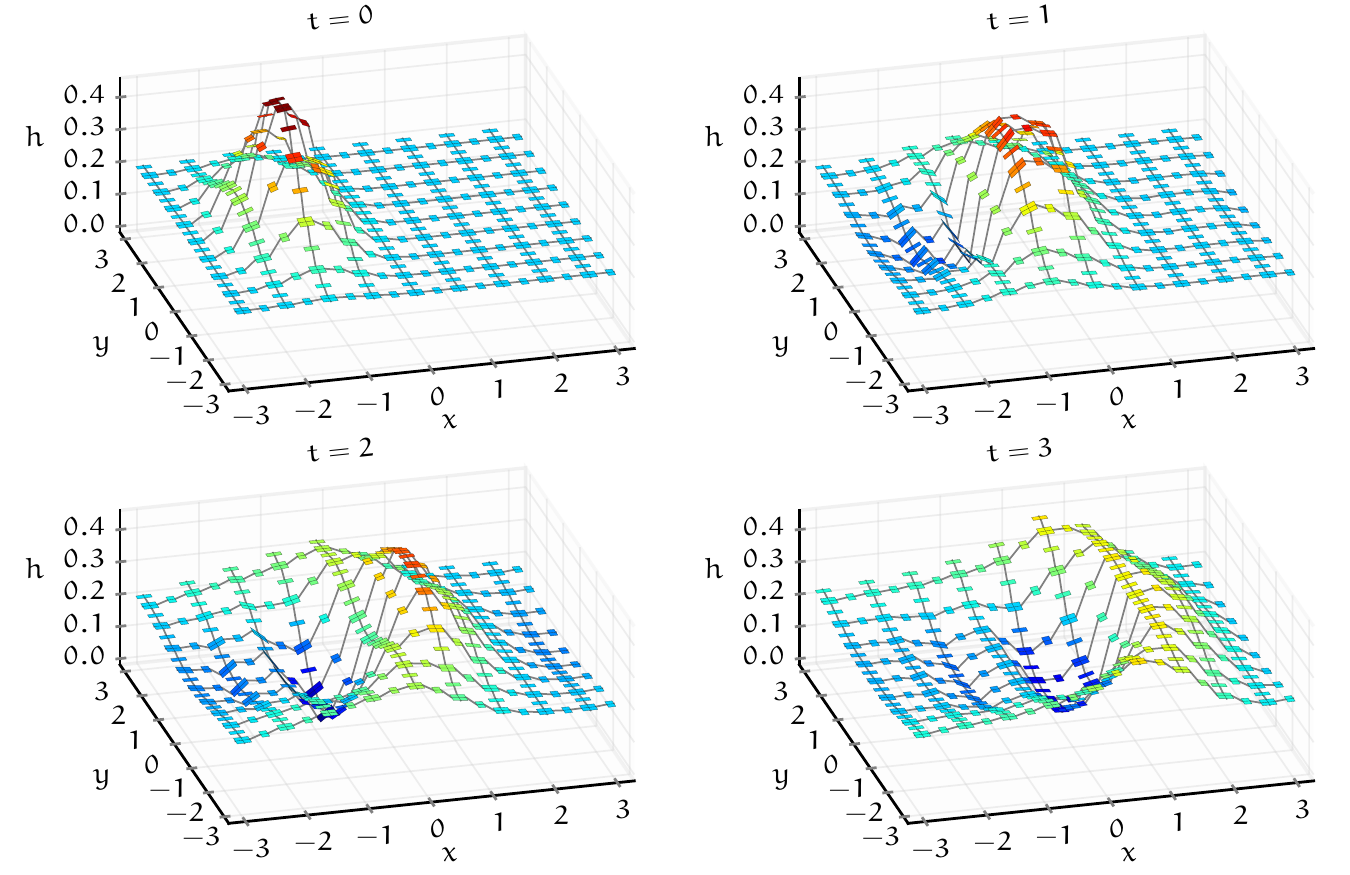}
\end{figure}

In the fluid dynamics of Earth's atmosphere and oceans, the length scales range from a few millimetres to several thousands of kilometres \parencite[p.~3]{Grooms2018_MltSclMdlsGeoInPhysclFluidDyn}.
We define a \emph{full-domain microscale simulation} as the detailed simulation (over all space-time scales) over the whole space-time simulation domain.
The main interest generally lies in large-scale dynamics only, yet the effect of the smallest scales that influence the emergent large-scale dynamics needs to be accounted for.
A full-domain microscale simulation over such a large space is impractical.
Hence, this article develops the foundation of an equation-free patch scheme that is practically feasible for such problems because it computes on only a small fraction of space (see \cref{fig:tSim2_gLnrWve_h_r0p1} for example). % [JD] Made "Equation-Free Patch Scheme" to "equation-free patch scheme" to make it consistent and because we are not introducing any acronym here.

Many multiscale modelling methods \parencite[e.g.,][]{Grooms2018_MltSclMdlsGeoInPhysclFluidDyn, Welsh2018_MltiSclMdlngFoodDryng, Emereuwa2020_MthmtclHmgnztnAndStchstcMdlngOfEnrgyStrgeSystms} aim to accurately model the macroscale physics by computing only within small coupled regions in the spatial domain.
The \emph{equation-free patch scheme} \cite[e.g.]{Kevrekidis2009_EqtnFreMltscleCmpttnAlgrthmsAndAplctns} is a flexible, computationally efficient, multiscale modelling approach.
Equation-free multiscale patch schemes have been developed, with proven consistency, and applied successfully, for dissipative systems \parencite{Roberts2005_HiOrdrAcrcyGapTthScme,Roberts2007_GnrlTthBndryCndnsEqnFreeMdlng,Bunder2017_GdCplingFrTheMltsclePtchSchmeOnSystmsWthMcrscleHtrgnty,Maclean2021_AtlbxOfEqtnFreFnctnsInMtlbOctveFrEfcntSystmLvlSmltn}.
Systems that predominantly describe waves pose significant challenges due to the wave dynamics being on the verge of instability.
\cite{Cao2013_MltSclMdlCplsWave, Cao2015_MltSclMdlCplsNonLrWave} extended the patch scheme to 1D wave-like systems using a staggered macroscale grid of patches in 1D space, where each patch itself contains a staggered microscale grid in 1D~space.
This article further develops the staggered patch scheme for waves in 2D~space in order to address and resolve crucial issues for multi-D wave simulation.
Detailed exploration over a wide range of parameters establishes stability, accuracy, consistency, computational efficiency, and insensitivity to numerical roundoff errors.

Crucially, patch schemes apply to, or wrap around, \emph{any given} detailed microscale computational function that a scientist or engineer trusts to encompass the multiscale physics of interest \cite[e.g.,][]{Maclean2021_AtlbxOfEqtnFreFnctnsInMtlbOctveFrEfcntSystmLvlSmltn}.
Patch schemes efficiently make macroscale predictions of a multiscale system by performing detailed microscale computation only within small sparsely distributed patches (e.g., \cref{fig:tSim2_gLnrWve_h_r0p1}).
Accurate macroscale predictions are made by correctly coupling the patches via appropriate interpolation over unsimulated macroscale space \parencite{Kevrekidis2004_EqFreeCompAiddAnlysMltiSclSys, Kevrekidis2009_EqtnFreMltscleCmpttnAlgrthmsAndAplctns, Hyman2005_PtchDynmcsForMltsclePrblms}.
There is no derived equation as a closure that describes a macroscale model, hence the name \emph{equation-free}. 
The macroscale information one obtains about the system is the computed data of the spatially distributed patches \parencite{Kevrekidis2009_EqtnFreMltscleCmpttnAlgrthmsAndAplctns}.
This article explores in detail two kinds of patch coupling which give two good families of equation-free patch schemes \parencites[\S3.1]{Divahar2022_AcrteMltiscleSmltnOfWveLkeSystms}[\S2.2.1, \S2.2.3]{Bunder2020_LrgeScleSmltnOfShlwWtrWvsVaCmpttnOnlyOnSmlStgrdPtchs}:
a \emph{spectral patch scheme} using global spectral interpolation (\cref{ssc:spctrlCplng});
and \emph{polynomial patch schemes} using local polynomial interpolation (\cref{ssc:plynmlCplng}).
One can achieve an arbitrarily high order of macroscale consistency for patch schemes, and hence controllable accuracy, via appropriate high-order interpolation for patch coupling \parencite{Roberts2005_HiOrdrAcrcyGapTthScme, Roberts2007_GnrlTthBndryCndnsEqnFreeMdlng}.

These multiscale patch schemes, by computing only within a small fraction of the whole domain, offer enormous computational savings in many physical applications.
%
% For example, the simulation of \cref{fig:tSim2_gLnrWve_h_r0p1} with its small patch ratio~\(r=0.01\) (ratio of patch width to patch separation) is 1300 times quicker than the corresponding simulation computed over the full spatial domain.
For example, a simulation as in \cref{fig:tSim2_gLnrWve_h_r0p1} but with a smaller patch ratio~\(r=0.01\) (ratio of patch width to patch separation) is 1300 times quicker than the corresponding (same resolution) simulation computed over the full spatial domain. % [JD] Figure~1 simulation is for r=0.1, not r=0.01, so changed to better clarify.
\cref{sec:cmptnlSvngs_gLnrWve} establishes scenarios for 2D waves where speed-ups of up to~\(10^5\) are achieved compared to a detailed full-domain simulation.
Thus, patch schemes have the potential to accurately and efficiently predict emergent macroscale waves in detail over large spatial domains from a given multiscale wave-like system (e.g., floods, tsunamis).

\begin{SCfigure}[50]
  \def\svgscale{0.85}
    \caption{\label{fig:flDmnGrd_stgrd}%
      Schematic microscale staggered grid where variables are stored only on staggered/alternating discrete points (\emph{nodes}~\(\hNode\,h\), \(\uNode\,u\), and~\(\vNode\,v\)).
      This staggered grid has \(6 \times 6\) (\(n=6\)) grid intervals and \(3 \times 3 = 9\) micro-cells ({\cpq orange} square).
      The unfilled nodes (\(\hNodeB\,h\), \(\uNodeB\,u\), and~\(\vNodeB\,v\)) indicate discrete \(n\)-periodic boundary values.
      An \(n \times n\) full-domain microscale staggered grid has \(n^2/4\)~micro-cells and \(3n^2/4\)~state variables.
    }
    \small
    \input{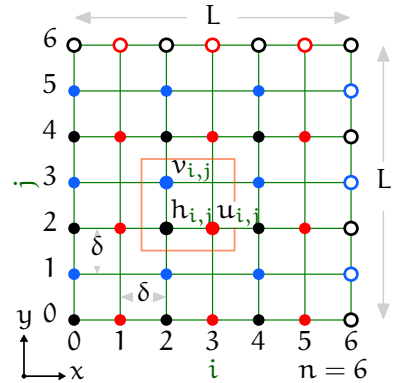}
\end{SCfigure}

Equation-free multiscale patch schemes have been developed, proven, and applied successfully for dissipative systems \parencite[e.g.,][]{Roberts2005_HiOrdrAcrcyGapTthScme, Roberts2007_GnrlTthBndryCndnsEqnFreeMdlng, Bunder2017_GdCplingFrTheMltsclePtchSchmeOnSystmsWthMcrscleHtrgnty, Maclean2021_AtlbxOfEqtnFreFnctnsInMtlbOctveFrEfcntSystmLvlSmltn}.
However, computational schemes for wave-like systems with small dissipation are often inaccurately unstable due to methodological quirks and/or roundoff errors (\cites[p.136]{Hinch2020_ThnkBfreYuCmpte_APrldeToCmpttnlFldDynmcs}[pp.~70--73]{Zikanov2010_EsntlCmpttnlFldDynmcs}[pp.~232--243]{Anderson1995_CmpttnlFldDynmcs}).
In order to represent the physical wave phenomena, a patch scheme for wave-like systems needs to navigate these issues---issues that are more difficult in multiple space dimensions.
For wave-like systems in full-domain modelling, a common strategy for accurate and robust spatial discretisation schemes is to use staggered spatial grids as shown schematically in \cref{fig:flDmnGrd_stgrd} \parencite[from Fig.~1]{Divahar2022_StgrdGrdsFrMltdmnsnlMltscleMdlng}.
Staggered grids preserve much of the wave characteristics \parencites[\S2]{Divahar2022_StgrdGrdsFrMltdmnsnlMltscleMdlng}
  [p.46, \S3.2]{Lauritzen2011_NmrclTchnqsFrGlblAtmsphrcMdls}
  [p.55, \S2.2.1]{Olafsson2021_UncrtntsInNmrclWthrPrdctn}
  [Figs.~8 \&~9]{Fornberg1999_SptlFnteDfrnceAprxmtnsFrWveTypeEqtns}
  {Fornberg1990_HghOrdrFnteDfrncsAndThePsdspctrlMthdOnStgrdGrds},
and typically support higher accuracy simulations compared to simulations of the same order on collocated grids.
Furthermore, the group velocity of the energy propagation in the numerical waves on a staggered grid is approximately in the correct direction, whereas on collocated grids \parencite[Fig.~1, left]{Divahar2022_StgrdGrdsFrMltdmnsnlMltscleMdlng} the group velocity for large wavenumbers is often in the opposite direction (\cites[p.46, \S3.2]{Lauritzen2011_NmrclTchnqsFrGlblAtmsphrcMdls}
[p.55, \S2.2.1]{Olafsson2021_UncrtntsInNmrclWthrPrdctn}).
To excellently preserve wave properties in 1D~space, \cite{Cao2013_MltSclMdlCplsWave} extended the patch scheme to a 1D staggered macroscale grid of patches, where each patch itself contains a 1D \text{staggered microscale grid.}

\begin{SCfigure}
    \caption{\label{fig:PtchGrd_n2t0}%
      A staggered patch grid is a \emph{macro-grid} of staggered patches ({\cIJ violet} squares) enclosing \emph{sub-patch micro-grids} ({\cij green}); {\cpq orange} squares containing three patches are \emph{macro-cells}.
      The distance separating patches (between the {\cIJ violet} lines) is the \emph{inter-patch spacing}~\(\Delta\).
      The sub-patch grids ({\cij green}) has \emph{micro-grid spacing}~\(\delta\).
      Each patch is boarded by two layers of \emph{edge nodes} whose values are determined by the scheme's inter-patch coupling.
    }
    \input{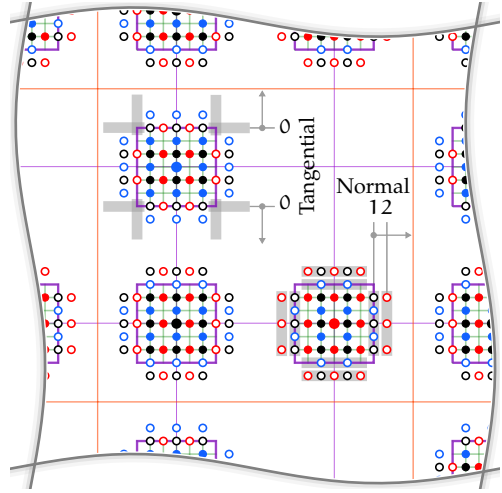}
\end{SCfigure}

\textcite{Divahar2022_StgrdGrdsFrMltdmnsnlMltscleMdlng} began extending the concept of staggered patch grids to multiple space dimensions by analysing all the possible~167\,040 2D staggered patch grids for wave-like systems.
Among the 167\,040 possible 2D multiscale staggered patch grids, \textcite[\S4.2, \S4.3]{Divahar2022_StgrdGrdsFrMltdmnsnlMltscleMdlng} showed that only 120~staggered patch designs constitute stable and accurate patch schemes for linear wave-like systems.
For analysing two families of equation-free multiscale patch schemes, this article focuses in detail on one of these~120, namely the staggered patch design depicted in \cref{fig:PtchGrd_n2t0} \parencite[patch grid~\#79985, Fig.~5]{Divahar2022_StgrdGrdsFrMltdmnsnlMltscleMdlng}.
Throughout the rest of this article, a \emph{staggered patch grid} refers to the specific multiscale staggered grid of \cref{fig:PtchGrd_n2t0}.

For wave-like systems in 2D~space, previous articles \parencite{Bunder2020_LrgeScleSmltnOfShlwWtrWvsVaCmpttnOnlyOnSmlStgrdPtchs, Divahar2022_StgrdGrdsFrMltdmnsnlMltscleMdlng} present promising preliminary results from the staggered patch schemes over the patch grid of \cref{fig:PtchGrd_n2t0}.
This article, for the generic case of weakly damped linear wave \pde{}s~\cref{eqs:PDEs_gLnrWve}, explores in great detail two families of the staggered patch schemes (\cref{secps2Dw}), for their 
accuracy and consistency (\cref{sec:ptchSchmsAreAcrte_gLnrWve,sec:ptchSchmsAreCnsstnt_gLnrWve}), 
insensitivity to numerical roundoff errors (\cref{sec:ptchSchmsAreNtSnstveToNmrclErrs_gLnrWve}), 
stability (\cref{sec:ptchSchmsAreStble_gLnrWve}), 
and the computational savings (\cref{sec:cmptnlSvngs_gLnrWve}).
Subsequent articles will explore the staggered patch schemes applied to the nonlinear wave \pde{}s for viscous and for turbulent `shallow' water flows.

% [JD] The following paragraph is essentially a merged version of what Trent suggested and Prof. Tony suggested earlier for my thesis (the paragraph just before Section 1.1).
An important objective of this article is to establish that the multiscale staggered patch schemes accurately and efficiently simulate \emph{wave-like systems} with small dissipation, and \emph{using a given microscale model}.
The focus on the weakly damped linear waves, despite no multiscale structure in their solution, is to understand the numerical characteristics of the patch schemes for a well-understood system that is amenable to analysis.
For the weakly damped linear waves, one could reasonably accurately use a coarse spatial grid, without a strong need for multiscale modelling.
But the aim is to use a given microscale code, and anticipate the patch schemes to efficiently simulate where fine micro-grids are essential (e.g., heterogeneities and turbulence).

\section{Patch schemes for 2D waves}
\label{secps2Dw}

In \cref{fig:PtchGrd_n2t0}, the ({\cIJ violet}) squares are small \emph{patches} enclosing ({\cij green}) \emph{sub-patch micro-grids}, whereas the ({\cpq orange}) squares containing a triangle of three patches are \emph{macro-cells}.
The side length of these square patches is the \emph{patch size}~\(l\).
The patches are placed on a ({\cIJ violet}) grid with \emph{inter-patch spacing}~\(\Delta\). 
Within every patch, the ({\cij green}) microscale grid has spacing~\(\delta\).

A finite-sized 2D patch grid over an \(L \times L\) spatial domain is specified by three parameters.
\begin{itemize}
  \item \(N\) is the number of macro-grid intervals ({\cIJ violet}) in the periodic domain in each of the~\(x\)- and \(y\)-directions.  
  Hence, the patch spacing~\(\Delta = L/N\)
  \item \(n:=l/\delta\) is the number of micro-grid intervals ({\cij green}) within a square patch in each of the~\(x\)- and \(y\)-directions. 
  \item The \emph{patch ratio}~\(r:=l/(2\Delta)\) quantifies the ratio of the simulated to the unsimulated space in each spatial dimension. 
  In practical use, patch ratios~\(r\) are small, typically ranging from~\(0.0001\) to~\(0.1\).
\end{itemize}
Throughout this article, we non-dimensionalise lengths in the problem with respect to the domain size~\(L\) so that non-dimensionally the spatial period is \(L = 2\pi\)\,. 
That is, herein we address solutions that are \(2\pi\)-periodic in space.
Hence, the patch spacing \(\Delta = L/N = 2\pi/N\), and the sub-patch micro-grid spacing~\(\delta = 2L r /(N n) = 4\pi r/(N n)\).\footnote{We use the same symbol~\(n\) and~\(\delta\) for both the full-domain micro-grid and the sub-patch micro-grid, and disambiguate by words and/or context.}

The unfilled circles in \cref{fig:PtchGrd_n2t0} are  patch \emph{edge nodes}.
For first-order \pde{}s, the patch schemes over this staggered patch grid interpolate field values to only those edge nodes on the four edges of the  ({\cIJ violet}) square patches in \cref{fig:PtchGrd_n2t0}. 
However, discretisation of \pde{}s with higher-order spatial derivatives (e.g., diffusion terms~\(\nabla^2 u, \nabla^2 v\) in \eqref{eqs:PDEs_gLnrWve}) or more complex terms (e.g., mixed derivatives) requires interpolating field values to additional layers of edge nodes just outside the ({\cIJ violet}) squares, as also plotted in \cref{fig:PtchGrd_n2t0} \parencite[\S3.3]{Divahar2022_StgrdGrdsFrMltdmnsnlMltscleMdlng}.

\subsection{Weakly damped, spatially discrete, wave equations}

We consider computational simulations arising from the generic non-dimensional \emph{weakly damped linear wave \pde{}s}
\begin{subequations} \label{eqs:PDEs_gLnrWve}
\begin{align}
  \label{eqn:PDE_gLnrWve_h}
  \doh{h}{t} &= - \doh{u}{x} - \doh v y,\\
  \label{eqn:PDE_gLnrWve_u}
  \doh{u}{t} &= - \doh{h}{x} - c_D u + c_V \doh[2]{u}{x} + c_V \doh[2]{u}{y} \,, \\
  \label{eqn:PDE_gLnrWve_v}
  \doh{v}{t} &= - \doh{h}{y} - c_D v \, + c_V \doh[2]{v}{x} \, + c_V \doh[2]{v}{y}\,,
\end{align}
\end{subequations}
with linear drag and viscous diffusion (viscosity) characterised by the \emph{dissipation coefficients}~\(c_D, c_V \geq 0\) respectively. 
The case \(c_D = c_V = 0\) corresponds to \emph{ideal waves} with no dissipation (\cites[pp.136--137]{Dean1991_WtrWveMchncsFrEngnrsAndScntsts}[pp.257--258]{Mehaute1976_AnIntro2HydrdynmcsWtrWvs}).
We focus on the computation of the dynamics, and so adopt the simple boundary conditions that the three fields~\(h\), \(u\), and~\(v\) are \(L\)-periodic in both~\(x\) and~\(y\).
The wave system~\cref{eqs:PDEs_gLnrWve} arises as a description of linear waves in many physical scenarios and so the patch schemes and results developed here apply very broadly.

We suppose that there are given microscale features to be resolved on the given length scale~\(\delta\).
In future research, \(\delta\)~will be the length scale required to resolve microscale heterogeneities and/or intricate sub-patch dynamics.
Consequently, this article assumes~\(\delta\) is a fixed given value.

Approximating the spatial derivatives of the weakly damped linear wave \pde{}s~\cref{eqs:PDEs_gLnrWve} by central finite differences on the \emph{staggered micro-grid} of \cref{fig:flDmnGrd_stgrd}, spacing~\(\delta\), gives the following microscale discretisation
\begin{subequations}  \label{eqs:FDEs_mN_gLnrWve}
\begin{align}
  { \color{dClr_hNode}\bullet\;} \frac{d h_{{ i,j}}}{d t} &=
    -\frac{u_{ i+1,j} - u_{i-1,j}}{2 \delta}
    -\frac{v_{i,j+1} - v_{i,j-1}}{2 \delta}
\qquad \text{ for \(i,j\) even;} 
  \label{eqn:FDEs_gLnrWve_h}
\\ 
  { \color{dClr_uNode}\bullet\;} \frac{d u_{{i,j}}}{d t} &=
    - \frac{
      h_{i + 1,j}
      - h_{i - 1,j}
      }{
      2 \delta
      }
    - c_D u_{{i,j}} 
    + c_V \frac{
        u_{i-2,j}
        - 2 u_{{i,j}}
        + u_{i+2,j}
        }{
        4 \delta^{2}
        }
        \notag\\&\quad{}
      + c_V \frac{
        u_{i,j-2}
        - 2 u_{{i,j}}
        + u_{i,j+2}
        }{
        4 \delta^{2}
        }
  \qquad \text{ for \(i\) odd, \(j\) even;}
  \label{eqn:FDEs_gLnrWve_u}
    \\
  { \color{dClr_vNode}\bullet\;}\, \frac{d v_{{i,j}}}{d t} &=
    - \frac{
        h_{i,j + 1}
        - h_{i,j - 1}
        }{
        2 \delta
        }
    - c_D v_{{i,j}} 
     + c_V
      \frac{
        v_{i-2,j}
        - 2 v_{{i,j}}
        + \, v_{i+2,j}
        }{
        4 \delta^{2}
        }
        \notag\\&\quad{}
      + c_V \frac{
        v_{i,j-2}
        - 2 v_{{i,j}}
        + v_{i,j+2}
        }{
        4 \delta^{2}
        }
    \qquad \text{ for \(i\) even, \(j\) odd.}
  \label{eqn:FDEs_gLnrWve_v}
\end{align}
\end{subequations}
The \emph{full-domain microscale model} is then to compute with~\eqref{eqs:FDEs_mN_gLnrWve} over the entire \(L \times L\) spatial domain.
A \emph{patch scheme} is to compute with~\eqref{eqs:FDEs_mN_gLnrWve} only in small, sparsely distributed, patches of space.
Crucially, the patch scheme is to apply to any computational model, such as the discretisation~\cref{eqs:FDEs_mN_gLnrWve}, not to the corresponding underlying \pde{}s~\cref{eqs:PDEs_gLnrWve}.
So all discussions of stability and accuracy of the patch scheme are relative to the discretisation~\cref{eqs:FDEs_mN_gLnrWve}, not the \pde{}s~\cref{eqs:PDEs_gLnrWve}.
The patch scheme adopts the view that the given computational model on some microscale~\(\delta\), such as the discretisation~\cref{eqs:FDEs_mN_gLnrWve}, is trusted to represent the multiscale physics of interest.

The \emph{staggered patch scheme} uses~\eqref{eqs:FDEs_mN_gLnrWve} on the \(n \times n\) sub-patch staggered micro-grids of the patches in~\cref{fig:PtchGrd_n2t0}.
The patches are arrayed on an \(N\times N\) macro-grid over the spatial domain.
The macro-grid is divided into \(N/2 \times N/2\) macro-cells with each ({\cpq orange})  macro-cell containing three patches: \(\hNode\,h\)-centred, \(\uNode\,u\)-centred, and \(\vNode\,v\)-centred arranged as in \cref{fig:PtchGrd_n2t0}.
The patch scheme is closed by defining a patch coupling that determines the edge values of each patch (open symbols in~\cref{fig:PtchGrd_n2t0}) from representative field values of surrounding patches.
To distinguish patch scheme quantities from full-domain quantities we use superscripts of the patch index~\(I,J\): for example, \(h^{I,J}_{i,j}\)~is the value of the \(h\)-field at the \((i,j)\)th~micro-grid point in the \((I,J)\)th~patch.

The patch scheme state vector~\(\xv^{I}\) is the collection of all the \emph{interior} node values (\(h^{I,J}_{i,j}\), \(u^{I,J}_{i,j}\), and~\(v^{I,J}_{i,j}\)), that is, not including the edge values. 
The patch edge values are determined by some interpolation function of the patch interior values, denoted as~\(\xv^{E}(\xv^{I})\)---the patch coupling. 
Then each staggered patch scheme may be represented as a dynamical system by a system of \ode{}s in the autonomous dynamical system form  
\begin{equation} \label{eqn:dynSys_pN_gnrcWve}
  \dby{\xv^I}{t} = \vec{F}\big( \xv^{I}; \xv^{E}(\xv^{I}) \big).
\end{equation}
Simulations by the patch scheme, such as \cref{fig:tSim2_gLnrWve_h_r0p1}, are obtained by numerically integrating the system~\eqref{eqn:dynSys_pN_gnrcWve}.
General characteristics of the patch scheme are determined by analysing the function~\(\vec F\) for various~\(\xv^{E}(\xv^{I})\): it is these general characteristics that we explore herein.

\Cref{ssc:spctrlCplng,ssc:plynmlCplng} detail two families of the staggered patch schemes based upon the patch coupling~\(\xv^{E}(\xv^{I})\): one family uses spectral interpolation (\cref{ssc:spctrlCplng}); and the other uses polynomial interpolation (\cref{ssc:plynmlCplng}).
The patch coupling~\(\xv^{E}(\xv^{I})\) computes patch edge values in two steps:
\begin{enumerate}
  \item  for each patch, from their respective interior microscale values, compute a \emph{macroscale patch value}---a representative aggregate value, also called amplitude or order parameter;
  \item  for each patch, compute its \emph{microscale edge values}  by interpolating from the macroscale values of neighbouring patches across the relatively large inter-patch spacing~\(\Delta\).
\end{enumerate}
Thus a chosen patch coupling~\(\xv^{E}(\xv^{I})\) provides the crucial two-way connection between the microscale and macroscale.

Each of the three types of patches in the staggered patch grid (\cref{fig:PtchGrd_n2t0}) has a centre node which is one of the three fields: the centre node is either an~\(\hNode\,h\), \(\uNode\,u\) or~\(\vNode\,v\) node.
Patch grids without a patch centre node for a field do not constitute stable staggered patch schemes \parencite[\S3.3]{Divahar2022_StgrdGrdsFrMltdmnsnlMltscleMdlng}, so are not considered herein.
We then define a patch's \emph{representative aggregate value}, denoted by~\(H_{I,J}\), \(U_{I,J}\) or~\(V_{I,J}\) respectively, as that of the patch's centre node.
We do \emph{not} invoke localised averages of sub-patch quantities, as is commonly done \cite[e.g.,][]{Carr2016, Liu2015, Kevrekidis09a}, since analysis proves that our centre-node definition can generally achieve arbitrarily high-order accuracy \cite[e.g.,][]{Roberts2011a, Bunder2020_LrgeScleSmltnOfShlwWtrWvsVaCmpttnOnlyOnSmlStgrdPtchs}.

Consequently, for an \(N\times N\) macro-grid of staggered patches, on an \(L\times L\) domain, we have three \(N/2 \times N/2\) arrays, one for each of~\(H\), \(U\), and~\(V\) macroscale `aggregate' values.
To complete the patch scheme, the outstanding issue is to determine the fields on the patch edges from these macroscale values: two choices are described by \cref{ssc:spctrlCplng,ssc:plynmlCplng}.

\subsection{Spectral patch scheme for best accuracy}
% --------------------------------------------------------------------
\label{ssc:spctrlCplng}

We are given a staggered patch grid (\cref{fig:PtchGrd_n2t0}) with three \(N/2 \times N/2\) arrays of \(H\), \(U\) and~\(U\) macroscale values from the three kinds of patches.
All patches of the same kind (\(\hNode\,h\), \(\uNode\,u\) or \(\vNode\,v\)-centred) are equally spaced with the \emph{inter-cell distance}~\(2\Delta\).
The \emph{spectral patch scheme} described here uses Fourier interpolation to compute the microscale patch edge values~\(\hNodeE\,h\), \(\uNodeE\,u\), and~\(\vNodeE\,v\) from the equispaced macroscale values~\(H\), \(U\), and~\(V\).

Consider the 2D discrete Fourier transform~(\dft) of the~\(N/2 \times N/2\) array~\(H\) of the macroscale \(h\)-field:
\begin{align} \label{eqn:DFT_h}
  \widetilde{H}_{k_x,k_y} = \dft(H)
    &:= \sum_{I,J=1}^{N/2} H_{I, J}\cdot \exp[ -\i (k_x x_I + k_y y_J) ],
\end{align}
where the wavenumbers \(k_x,k_y \in \{-(N/2 - 1)/2,\ldots,(N/2 - 1)/2\}\), and similarly for the \(U,V\)~arrays.
The 2D inverse semidiscrete Fourier transform of the discrete Fourier transform~\(\widetilde{H}\), an \(N/2 \times N/2\) array, gives a continuous function \(h(x,y)\) which is the interpolated macroscale field at  arbitrary position~\((x,y)\):
\begin{align*} %\label{eqn:ISDFT_hfield}
  h(x,y) = \isdft(\widetilde{H})
    &:= \frac{1}{(N/2)^2}\sum_{k_x,k_y} \widetilde{H}_{k_x,k_y}  \cdot \exp[ \i (k_x x + k_y y) ].
\end{align*}
To compute the interpolated \(h\)-field at a position~\((\xi,\eta) := \big[(x,y)-(x_I,y_J)\big] /\Delta\) relative to the patch centres~\((x_I,y_J)\) and scaled by the inter-patch spacing~\(\Delta\), for every patch~\((I,J)\), we thus compute the shifted-inverse 
\begin{align} \label{eqn:ISDFT_h}
  h^{I,J}(\xi,\eta) 
    &= \frac{1}{(N/2)^2}\sum_{k_x,k_y} \widetilde{H}_{k_x,k_y}\e^{\i(k_x\xi+k_y\eta)\Delta} \cdot \exp[ \i (k_x x_I + k_y y_I) ]
\end{align}
via the inverse semidiscrete Fourier transform of the array~\(\widetilde{H}_{k_x,k_y}\e^{\i(k_x\xi+k_y\eta)\Delta}\).
Similarly for the \(u,v\)-fields.
When \(N/2\) is even, special handling is required for the Nyquist highest frequency component.  We avoid such special handling by requiring \(N/2\) to be odd (\(N \in \{6, 10, 14, 18, \ldots\}\)).

The Fast Fourier Transform~(\fft) provides efficient calculation of the transforms~\cref{eqn:DFT_h,eqn:ISDFT_h}.
We need three~\fft{}s for~\cref{eqn:DFT_h}, one for each field, and the number of inverse-\fft{}s~\eqref{eqn:ISDFT_h} is three times the number of edge nodes around one patch.
By capturing the global information from all macroscale wave components, spectral interpolation achieves high accuracy \parencite[e.g.,][]{Roberts2011a, Bunder2020_LrgeScleSmltnOfShlwWtrWvsVaCmpttnOnlyOnSmlStgrdPtchs}.

Such spectral interpolation is typically restricted to periodic macroscale boundary conditions on rectangular domains, but the cognate Chebyshev interpolation may promise high accuracy with more general boundary conditions and domains.

\subsection{Polynomial patch schemes for complex geometry}
% --------------------------------------------------------------------
\label{ssc:plynmlCplng}

For general domain shapes and general boundary conditions, local polynomial interpolation is generally more widely applicable than the spectral interpolation of \cref{ssc:spctrlCplng}.
Consequently, this section develops a family of \emph{polynomial staggered patch schemes}, named {square-p\(p\)} for integer~\(p\), whose patch coupling is 2D Lagrangian polynomial interpolation over a near-square region.
The parameter~\(p\) denotes the order of the interpolating polynomial.

\begin{figure}
\caption{\label{fig:egIntrpltnStncls}%
  Example interpolation stencils of the polynomial patch scheme, of low order~\(p\), for interpolating \(\vNodeE\,v\) of \(\hNode\,h\)-, \(\uNode\,u\)-, \(\vNode\,v\)-centred patches in \cref{fig:PtchGrd_n2t0}.
  }
\centering
\begin{subfigure}{0.31\textwidth}
  \centering
  \caption{\label{fig:hEdgeNdeIntrpStncl-square-p2}{square-p2}, \(\hNode\,h\)-patch}
  \vspace{-1ex}
  \includegraphics[scale=1, trim={3.69in, 3.5in, 0, 11}, clip]{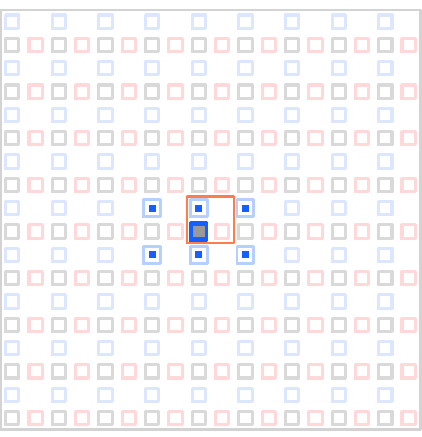}
\end{subfigure}
~
\begin{subfigure}{0.31\textwidth}
  \centering
  \caption{\label{fig:uEdgeNdeIntrpStncl-square-p2}{square-p2}, \(\uNode\,u\)-patch}
  \vspace{-1ex}
  \includegraphics[scale=1, trim={3.69in, 1.75in, 0, 137}, clip]{figs/interpStencilPlots/square-p2}
\end{subfigure}
~
\begin{subfigure}{0.31\textwidth}
  \centering
  \caption{\label{fig:vEdgeNdeIntrpStncl-square-p2}{square-p2}, \(\vNode\,v\)-patch}
  \vspace{-1ex}
  \includegraphics[scale=1, trim={3.69in, 0in, 0, 263}, clip]{figs/interpStencilPlots/square-p2}
\end{subfigure}
\\[1ex]
\begin{subfigure}{0.31\textwidth}
  \centering
  \caption{\label{fig:hEdgeNdeIntrpStncl-square-p4}{square-p4}, \(\hNode\,h\)-patch}
  \vspace{-1ex}
  \includegraphics[scale=1, trim={3.69in, 3.5in, 0, 11}, clip]{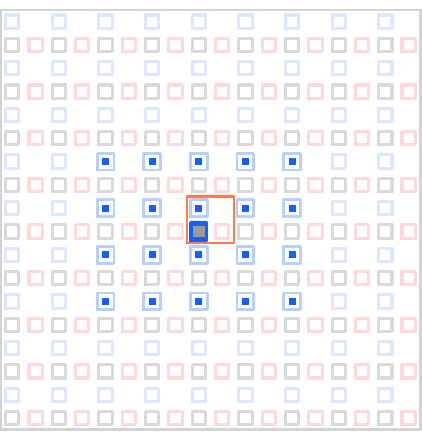}
\end{subfigure}
~
\begin{subfigure}{0.31\textwidth}
  \centering
  \caption{\label{fig:uEdgeNdeIntrpStncl-square-p4}{square-p4}, \(\uNode\,u\)-patch}
  \vspace{-1ex}
  \includegraphics[scale=1, trim={3.69in, 1.75in, 0, 137}, clip]{figs/interpStencilPlots/square-p4}
\end{subfigure}
~
\begin{subfigure}{0.31\textwidth}
  \centering
  \caption{\label{fig:vEdgeNdeIntrpStncl-square-p4}{square-p4}, \(\vNode\,v\)-patch}
  \vspace{-1ex}
  \includegraphics[scale=1, trim={3.69in, 0in, 0, 263}, clip]{figs/interpStencilPlots/square-p4}
\end{subfigure}
\end{figure}

Local polynomial interpolation computes the edge values of each patch using macroscale values of only some neighbouring patches.
The neighbourhood of a patch, characterised by the interpolation stencil, is of different sizes leading to different order~\(p\) of the interpolating polynomial.
For example, \cref{fig:egIntrpltnStncls} illustrates some interpolation stencils for low orders~\(p\) of the polynomial staggered patch schemes.
\textcite[Figs.~2.2.2--2.2.5]{Divahar2022_AcrteMltiscleSmltnOfWveLkeSystms} presents all the interpolation stencils (for coupling \(\hNodeE\,h\), \(\uNodeE\,u\), \(\vNodeE\,v\) edge values of \(\hNode\,h\)-, \(\uNode\,u\)-, \(\vNode\,v\)-centred patches) for each order \(p=2,4,6,8\).
All in the family of polynomial staggered patch schemes have stencils that are square or near-square (\cref{fig:egIntrpltnStncls}), but differ in size.
But depending upon the type of edge nodes (\(\hNodeE\,h\), \(\uNodeE\,u\), \(\vNodeE\,v\)) being interpolated and the type of the patch (\(\hNode\,h\)-, \(\uNode\,u\)-, \(\vNode\,v\)-centred) for which edge nodes are interpolated, some of the stencils are near-square rectangles.
For example, all the leftmost stencils in \cref{fig:egIntrpltnStncls} are near-square.

We define the \emph{polynomial interpolation order}~\(p\) as the maximum degree of the independent variables in the 2D Lagrangian basis polynomials of all the interpolation stencils of a staggered patch scheme.
For interpolation to the \((I,J)\)th~patch, we perform standard bivariate Lagrange interpolation (e.g., \cites%
  [\S10.10]{Gupta2019_NmrclMthds_FndmntlsAndAplctns}%
  [\S3.6]{Jain2004_NmrclMthds_PrblmsAndSltns}%
  [\S10.1]{Fletcher2020_SmiLgrngnAdvctnMthdsAndThrAplctnsInGscnce}%
  )
in patch local coordinate~\((\xi,\eta) := \big[(x,y)-(x_I,y_J)\big] /\Delta\).

\begin{table}
  \renewcommand{\arraystretch}{1.2}
  \caption{\label{tbl:ipPolyTble_v_v_square-p2}%
    2D Lagrangian basis polynomials \(\mathcal{B}_S(\xi, \eta)\) for the example of square-p2 patch coupling, and for interpolating the \(\vNode\,v\)-field from \(\vNode\,v\)-centred patches to the \(\vNodeE\,v\)~edge values of the middle patch, \(S=4\)\,.
    The stencil index~\(S \in \{0, 1, \ldots, n_S-1{=}8\}\), and the patch local coordinates~\((\xi,\eta):=\big[(x,y)-(x_I,y_J)\big]/\Delta\).
  }
\begin{equation*}
    \raisebox{-\height/2}{%
    \includegraphics[scale=0.9]{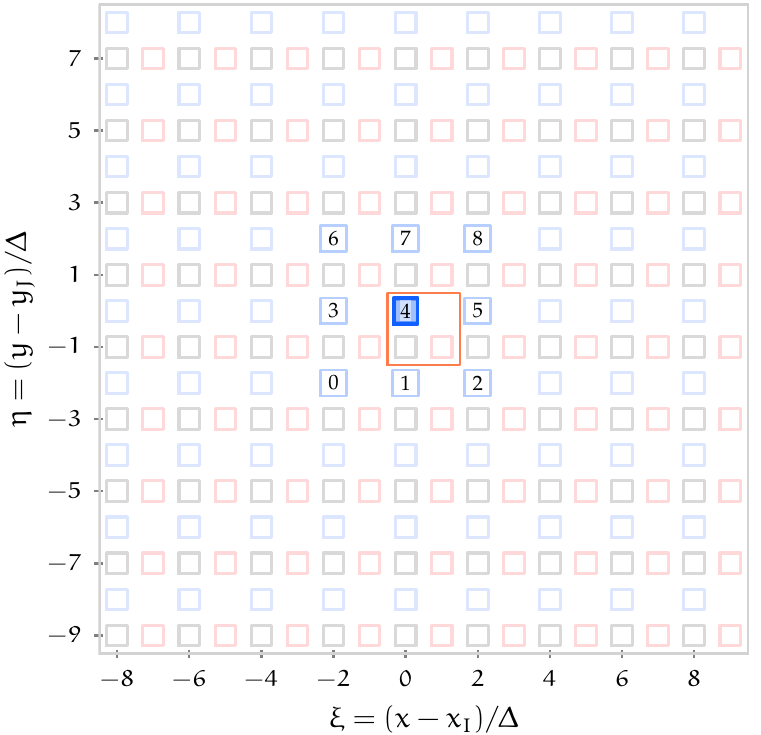}
    }
\quad
    \begin{array}{rl}
    \hline
    S
      & \multicolumn{1}{c}{\mathcal{B}_S(\xi, \eta)} \\ \hline
    0
      & \frac{\xi - 2}{4} \cdot \frac{\xi}{2} \cdot \frac{\eta - 2}{4} \cdot \frac{\eta}{2} \\
    1
      & - \frac{\xi - 2}{2} \cdot \frac{\xi + 2}{2} \cdot \frac{\eta - 2}{4} \cdot \frac{\eta}{2} \\
    2
      & \frac{\xi}{2} \cdot \frac{\xi + 2}{4} \cdot \frac{\eta - 2}{4} \cdot \frac{\eta}{2} \\
    3
      & - \frac{\xi - 2}{4} \cdot \frac{\xi}{2} \cdot \frac{\eta - 2}{2} \cdot \frac{\eta + 2}{2} \\
    4
      & \frac{\xi - 2}{2} \cdot \frac{\xi + 2}{2} \cdot \frac{\eta - 2}{2} \cdot \frac{\eta + 2}{2} \\
    5
      & - \frac{\xi}{2} \cdot \frac{\xi + 2}{4} \cdot \frac{\eta - 2}{2} \cdot \frac{\eta + 2}{2} \\
    6
      & \frac{\xi - 2}{4} \cdot \frac{\xi}{2} \cdot \frac{\eta}{2} \cdot \frac{\eta + 2}{4} \\
    7
      & - \frac{\xi - 2}{2} \cdot \frac{\xi + 2}{2} \cdot \frac{\eta}{2} \cdot \frac{\eta + 2}{4} \\
    8
      & \frac{\xi}{2} \cdot \frac{\xi + 2}{4} \cdot \frac{\eta}{2} \cdot \frac{\eta + 2}{4} \\
    \hline
    \end{array}
  \end{equation*}
\end{table}%
The bivariate Lagrange interpolation polynomial used for the patch coupling of all the polynomial staggered patch schemes is computed as
  \(f(\xi,\eta) = \sum_{S=0}^{n_S-1} \mathcal{B}_S(\xi, \eta) f_S\)\,,
where \(S\)~indexes the neighbouring patches in the chosen stencil, and
in terms of standard 2D basis polynomials~\(\mathcal{B}_S\), such as the example of \cref{tbl:ipPolyTble_v_v_square-p2}, and the known patch macroscale values~\(f_S\).

We restrict attention to near square macroscale stencils, such as \cref{fig:egIntrpltnStncls}.
Other non-`square', stencils may be worth investigating in the future.

\section{Staggered patches accurately resolve macroscale waves}
% ++++++++++++++++++++++++++++++++++++++++++++++++++++++++++++++++++++
\label{sec:ptchSchmsAreAcrte_gLnrWve}

This section shows that the spectral (\cref{ssc:spctrlCplng}) and the family of polynomial (\cref{ssc:plynmlCplng}) staggered patch schemes are accurate for the generic, weakly damped, linear wave system~\cref{eqs:FDEs_mN_gLnrWve}.
We explore a range of orders of interpolation~\(p\),  parameters of the patch grid~(\(N\), \(n\)), and parameters of the physical system~(\(C_D\), \(C_V\)).
We show the accuracy of the patch schemes for all simulations (as opposed to just computing for a few initial conditions) by comparing their eigenvalues with the eigenvalues of the `given' microscale model~\cref{eqs:FDEs_mN_gLnrWve} over the full domain (\cref{fig:flDmnGrd_stgrd}).

The eigenvalues are sufficient for measuring accuracy due to the following argument in the case of periodic boundary conditions.
Firstly, recall that for a linear, autonomous dynamical system, such as the system~\eqref{eqs:FDEs_mN_gLnrWve} considered herein, the general solution of the system (i.e., patch interior values~\(\xv^I\)) is generically a linear combination of the form \(\xv^I=\sum_k c_k\vec q_k\e^{\lambda_kt}\) for a complete set of eigenvectors~\(\vec q_k\) and eigenvalues~\(\lambda_k\), and for arbitrary constants~\(c_k\).
The eigenvalues and eigenvectors are of the linear operator (the Jacobian) encoding both the microscale details, such as~\eqref{eqs:FDEs_mN_gLnrWve}, and the specific patch coupling (\cref{ssc:spctrlCplng,ssc:plynmlCplng}).
Since the system is linear, we compute the Jacobian by concatenating the columns of the response of the patch scheme to unit impulse states \parencite[\S3.2.5]{Divahar2022_AcrteMltiscleSmltnOfWveLkeSystms}.
Secondly, our multiscale patch schemes (\cref{secps2Dw}) are translationally invariant to macroscale shifts of~\(2\Delta\) in both \(x,y\)-directions.
Because of this translational invariance, all the eigenvectors~\(\vec q_k\) have a sinusoidal modulation across the distributed macroscale patches.
These macroscale sinusoids exactly match those of the large-scale wave eigenvectors in the given full-domain system~\eqref{eqs:FDEs_mN_gLnrWve}.
Hence, in the case of a domain with periodic boundary conditions, the only error in the \emph{macroscale} predictions of the patch scheme are in the eigenvalues.

\subsection{Physical interpretation of spectra}

\begin{figure}
\centering
\caption{\label{fig:eigPlt_gLnrWve_square-p8_N10_n6_r0p1_cD1e-06_cV0p0001}%
typical structure of the spectrum of eigenvalues of the staggered patch scheme for weakly dissipative waves (here \(N = 10\), \(n = 6\), \(r=0.1\), \(\delta = \pi/150\), \(c_D=10^{-6}\), \(c_V=10^{-4}\)).  The spectrum typically has seven clusters: the right-hand four clusters are (75) macroscale modes; the left-hand three clusters are (1400) microscale sub-patch modes.  All the complex plane plots in this article utilise a quasi-log nonlinear scale, via \(\arcsinh\), to reasonably display the multiscale range of eigenvalues.}
\includegraphics{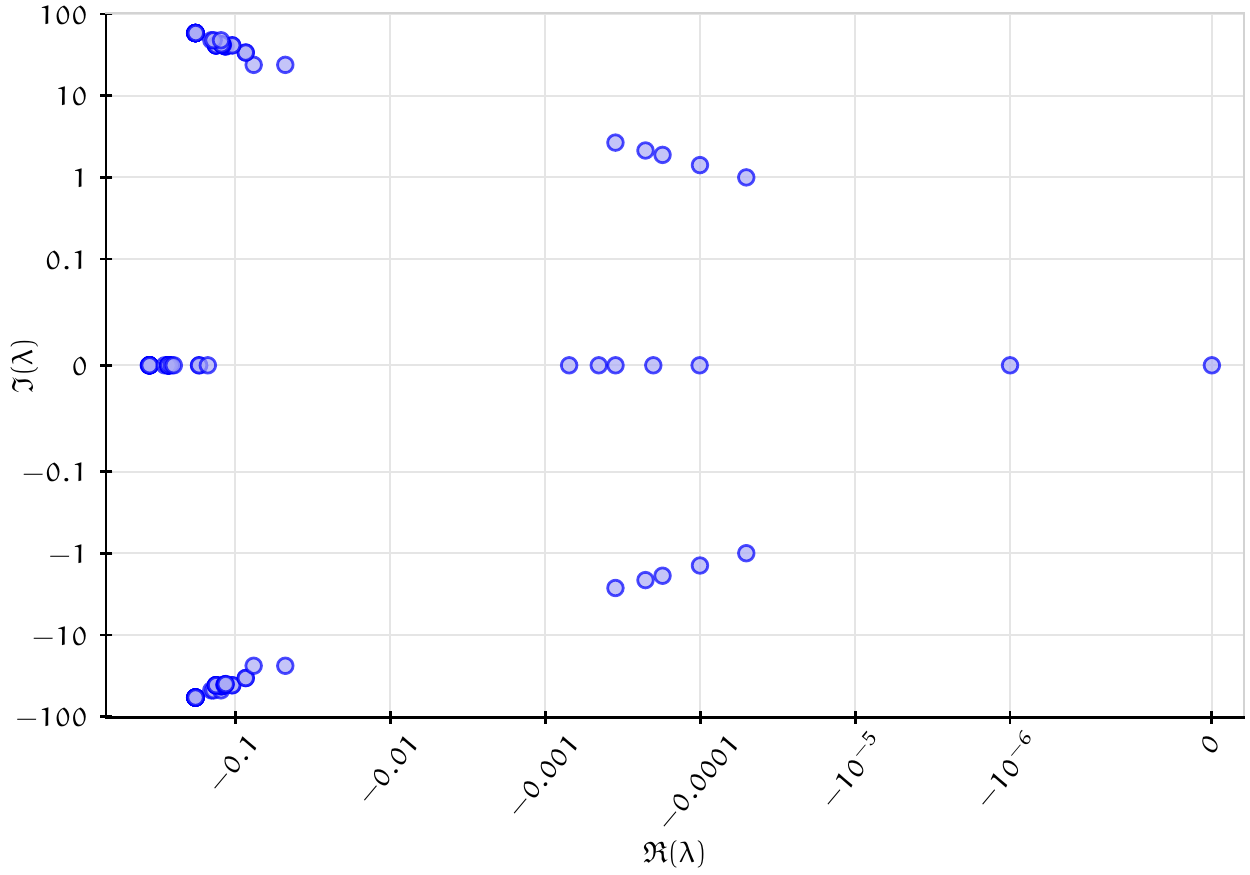} 
\end{figure}

Here we explain the structure of the spectra of patch scheme eigenvalues---introduced by the example plotted in \cref{fig:eigPlt_gLnrWve_square-p8_N10_n6_r0p1_cD1e-06_cV0p0001}. 
From such eigenvalues, we subsequently proceed to discuss the patch schemes' accuracy and their dependence upon the design parameters.

\needspace{2\baselineskip}
The dynamics of the patch scheme contains identifiable modes at both micro- and macro-scales.
\begin{itemize}
  \item \emph{Macroscale modes} are those eigenvectors that have macroscale spatial variation with little microscale structure within each patch.  Consequently, the macroscale modes have relatively small wavenumber.
  The corresponding eigenvalues are termed \emph{macroscale eigenvalues} (e.g., eigenvalues with \(\Re(\lambda) > -0.001\) in \cref{fig:eigPlt_gLnrWve_square-p8_N10_n6_r0p1_cD1e-06_cV0p0001}).
  \item \emph{Microscale modes} are those eigenvectors that have significant microscale structure irrespective of whether it is modulated by some macroscale variation.
  The corresponding eigenvalues are termed \emph{microscale eigenvalues} (e.g., eigenvalues with \(\Re(\lambda) < -0.01\) in \cref{fig:eigPlt_gLnrWve_square-p8_N10_n6_r0p1_cD1e-06_cV0p0001}).
\end{itemize}
Plots of various eigenvectors by \textcite[\S3.2.6]{Divahar2022_AcrteMltiscleSmltnOfWveLkeSystms} illustrate the above classification into microscale and macroscale modes.
Recall that we define whether an eigenvalue is macroscale or microscale by the spatial structure of its eigenvector, not by the properties of the eigenvalue itself.
However, the patch scheme modes are usually easily distinguishable by the eigenvalues alone as they form physically interpretable clusters (e.g., \cref{fig:eigPlt_gLnrWve_square-p8_N10_n6_r0p1_cD1e-06_cV0p0001,fig:compEigs_Sp_cDcV}).

The (\(3N^2/4\)) macroscale modes form the following four, physically interpreted, broad clusters of eigenvalues (\(\Re(\lambda) > -0.001\) in \cref{fig:eigPlt_gLnrWve_square-p8_N10_n6_r0p1_cD1e-06_cV0p0001}, and similarly in subsequent spectra).
\begin{itemize}
  \item A cluster of three near-zero eigenvalues (the two dots in \cref{fig:eigPlt_gLnrWve_square-p8_N10_n6_r0p1_cD1e-06_cV0p0001} with \(\Re(\lambda) > -10^{-5}\)) consists of one \(\lambda=0\) representing conservation of water~\(h\) (\(u = v = 0\)), and a pair of  eigenvalues \(\lambda = -c_D\) representing uniform flow in 2D space (\(h=0\)) decaying slowly due to drag.
  \item Two complex conjugate clusters (\(0.99<|\Im(\lambda)|<3\) in \cref{fig:eigPlt_gLnrWve_square-p8_N10_n6_r0p1_cD1e-06_cV0p0001}) represent macroscale waves---waves damped weakly by drag and viscosity.
   They arise from the \(N^2/4 - 1\) pairs of complex conjugate eigenvalues representing waves with non-zero wavenumbers \(|k_x|,|k_y|<N\pi/L\).
  \item A cluster of \(N^2/4 - 1\) real eigenvalues (\(-0.001<|\lambda|\leq-0.0001\) in \cref{fig:eigPlt_gLnrWve_square-p8_N10_n6_r0p1_cD1e-06_cV0p0001}) represents macroscale vortices in the \(xy\)-plane (\(h=0\))---weakly damped by drag and viscosity.  
  Such vortices may combine to form arbitrary, large-scale, \(xy\)-circulations.
\end{itemize}

The microscale modes generally form three broad clusters of eigenvalues  (\(\Re(\lambda) < -0.01\) in \cref{fig:eigPlt_gLnrWve_square-p8_N10_n6_r0p1_cD1e-06_cV0p0001}, and similarly in subsequent spectra).
The relatively large decay rate of all of these modes implies the patch system relatively quickly settles onto a slow manifold of the macroscale modes (\cites
  [\S5.3, p.302]{Roberts2003_LwDmnsnlMdlngOfDynmclSystmsApldToSmeDsptveFldMchncs}
  []{Zagaris2009_AnlyssOfTheAcrcyAndCnvrgnceOfEqtnFrePrjctnToASlwMnfld}
  []{Roberts1988_TheAplctnOfCntreMnfldThryToTheEvltnOfSystmWhchVrySlwlyInspce}
  []{Lorenz1986_OnTheExstnceOfASlwMnfld}).
\begin{itemize}
  \item Two complex conjugate clusters (\(10 <|\Im(\lambda)| < 100\) in \cref{fig:eigPlt_gLnrWve_square-p8_N10_n6_r0p1_cD1e-06_cV0p0001}), each containing \((N^2/4)(3n^2/4 - n - 1)\) eigenvalues, are of sub-patch microscale waves of large wavenumber.  
  These sub-patch microscale modes, with wave energy `bouncing around' within a patch while `leaking' to neighbouring patches, are relatively rapidly dissipated by viscosity.  

  \item A cluster of \((N^2/4)(3n^2/4 - 2n + 1)\) real eigenvalues (\(-1<\lambda<-0.1\) in \cref{fig:eigPlt_gLnrWve_square-p8_N10_n6_r0p1_cD1e-06_cV0p0001}) represents sub-patch microscale vortices (\(h=0\)), weakly linked to neighbouring patches, and rapidly dissipated by viscosity.  

\end{itemize}

\subsection{Compare with corresponding microscale model}

Recall that the objective of the multiscale patch schemes is to make \emph{accurate macroscale predictions} for some given microscale computational model.
Hence we compare the macroscale eigenvalues~\(\lambda_{p,M}\) of a patch scheme with the corresponding eigenvalues~\(\lambda_\mu\) of the microscale model~\eqref{eqs:FDEs_mN_gLnrWve} over the full-domain.
Hereafter, the subscripts
  \((\argdot)_p\)~refer to a patch scheme,
  \((\argdot)_{p,M}\)~refer to its macroscale modes,
  \((\argdot)_{p,\mu}\)~refer to its microscale modes, and
  \((\argdot)_\mu\)~refer to the full-domain microscale model (computed analytically), and
when necessary, superscripts \((\argdot)^A,(\argdot)^N\) distinguish analytically and numerically computed eigenvalues.
For accuracy considerations, the sub-patch microscale eigenvalues~\(\lambda_{p,\mu}\) are not of interest, except that they must not be unstable and should be separated from the slow macroscale eigenvalues.

\begin{figure}
  \caption{\label{fig:compEigs_Sp_cDcV}%
    \emph{Spectral patch scheme is exact} (\(N = 10\), \(n = 6\), \(r=0.1\), \(\delta = \pi/150\)), macroscale eigenvalues~\(\lambda_{p,M}\) agree with corresponding eigenvalues~\(\lambda_\mu\) of full-domain microscale model visually exactly (quantitatively within~\(10^{-12}\)).
    }
  \centering
  \begin{subfigure}[b]{0.489\textwidth}
    \centering
    \caption{\label{fig:compEigs_Sp_cD0cV0}%
      \(c_D = 0\), \(c_V = 0\).
      }
    \includegraphics[scale=0.95]{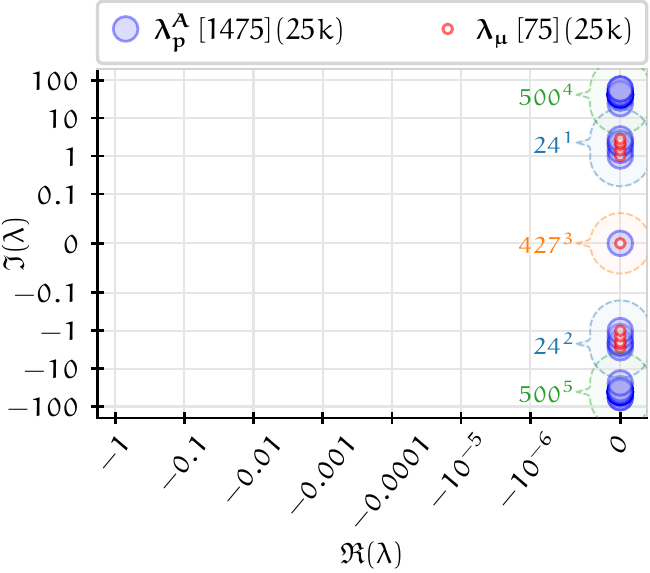}
  \end{subfigure}
  \begin{subfigure}[b]{0.489\textwidth}
    \centering
    \caption{\label{fig:compEigs_Sp_cD1e-6cV0}%
      \(c_D = 10^{-6}\), \(c_V = 0\).
      }
    \includegraphics[scale=0.95]{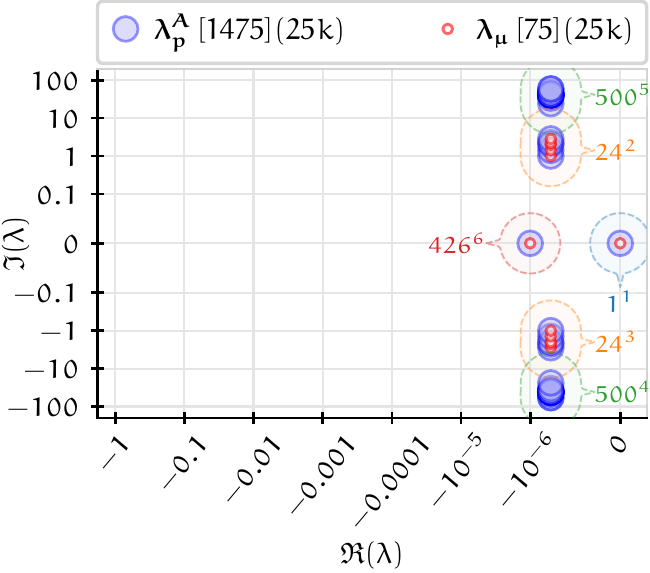}
  \end{subfigure}
  \begin{subfigure}[b]{0.489\textwidth}
    \centering
    \caption{\label{fig:compEigs_Sp_cD0cV1e-4}%
      \(c_D = 0\), \(c_V = 10^{-4}\).
      }
    \includegraphics[scale=0.95]{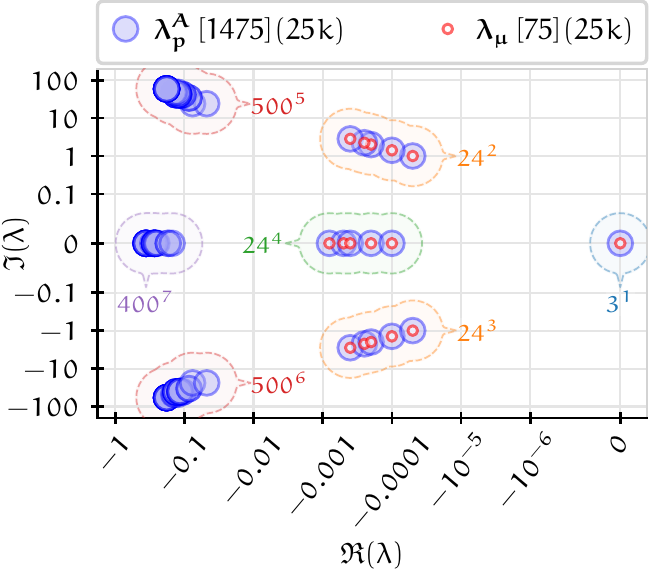}
  \end{subfigure}
  \begin{subfigure}[b]{0.489\textwidth}
    \centering
    \caption{\label{fig:compEigs_Sp_cD1e-6cV1e-4}%
      \(c_D = 10^{-6}\), \(c_V = 10^{-4}\).
      }
    \includegraphics[scale=0.95]{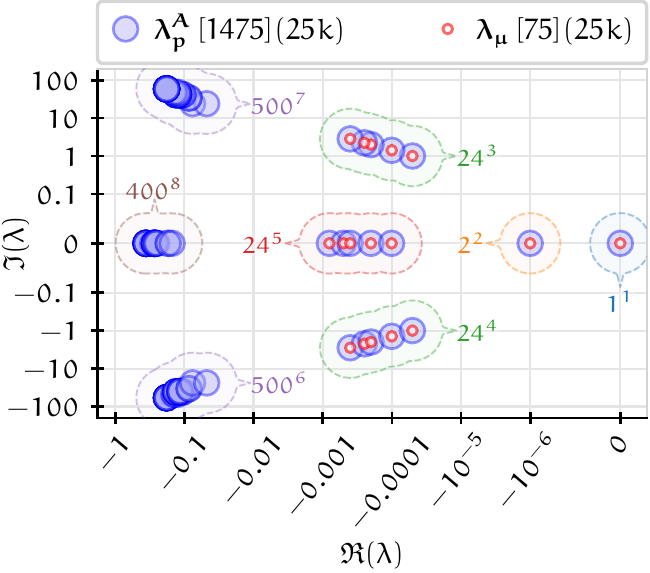}
  \end{subfigure}
\end{figure}

For the microscale discrete system~\eqref{eqs:FDEs_mN_gLnrWve}, \textcite[\S4.1]{Divahar2022_StgrdGrdsFrMltdmnsnlMltscleMdlng} derived the full-domain microscale eigenvalues
\begin{align}& \label{eqn:eig_muA}
\lambda_\mu  =
\begin{cases}
  - c_D + c_V \omega_{\mu0}^2\,,\\
  - \left( c_D + c_V \omega_{\mu0}^2 \right)/2
    \pm \i \sqrt{
      \omega_{\mu0}^2
      - \left[ \left( c_D + c_V \omega_{\mu0}^2 \right)/2 \right]^2
      } \,,
\end{cases}
\end{align}
where \(\omega_{\mu0} := \sqrt{ \sin^2{\left(k_x \delta \right) }/\delta^2 + \sin^2{\left(k_y \delta \right) }/\delta^2 } \) is the frequency of the discrete ideal wave (\(c_D = c_V = 0\)).
For each wavenumber~\((k_x,k_y)\), \cref{eqn:eig_muA}~gives three eigenvalues, one real and one complex conjugate pair, physically representing a vortex mode and a wave mode, respectively.
For subsequent complex eigenvalue spectra (\crefrange {fig:compEigs_Sp_cDcV} {fig:compEigs_Sq4_Nn}) we plot eigenvalues~\(\lambda_\mu\) (small red circles) computed by~\cref{eqn:eig_muA}, for all the \(3N^2/4\) \emph{macroscale} wavenumbers resolved by the discussed staggered patch scheme.
These are the \emph{reference eigenvalues} to assess the accuracy of a patch scheme's macroscale modes.

The eigenvalue plots \crefrange{fig:compEigs_Sp_cDcV}{fig:compEigs_Sq4_Nn} visually compare a patch scheme's eigenvalues~\(\lambda_p\) for~\cref{eqs:PDEs_gLnrWve} with the eigenvalues~\(\lambda_\mu\) of the full-domain microscale model~\cref{eqs:PDEs_gLnrWve}.
In the complex plane plots, the number within bracket~\([\argdot]\) in the legends gives the total number of eigenvalues for each set.
The number within parenthesis~\((\argdot k)\) in the legends gives the total number of wavenumbers for which the analytical eigenvalues~\(\lambda^A\) are evaluated (the~\(k\) here denotes the wavenumber~\(k\), not kilo-).

\subsection{Spectral coupling is highly accurate}
\label{ssc:SpctrlCplngIsHghlyAcrte}

\cref{fig:compEigs_Sp_cDcV} suggests that the spectral patch scheme is exact for the macroscale dynamics.
\Cref{fig:compEigs_Sp_cDcV} shows that, for all four combinations of \(c_D \in \{0, 10^{-6}\}\) and \(c_V \in \{0, 10^{-4}\}\),  the plotted full-domain eigenvalues~\(\lambda_\mu\) (red circles) visually exactly match each of the macroscale patch eigenvalues~\(\lambda_{p,M}\) (blue discs) of the spectral patch scheme.
Numerically we find such eigenvalues match to within~\(10^{-12}\) (i.e., to round-off error).
That is, the \emph{spectral staggered patch scheme makes effectively exact predictions of the macroscale dynamics}.

\cref{fig:compEigs_Sp_cD0cV1e-4,fig:compEigs_Sp_cD0cV1e-4}, with non-zero viscosity, shows the seven previously discussed clusters (although \cref{fig:compEigs_Sp_cD0cV1e-4} splits the three \(\lambda\approx 0\) into two `clusters').
However, in the cases of zero viscosity \cref{fig:compEigs_Sp_cD0cV0,fig:compEigs_Sp_cD1e-6cV0} show that the eigenvalue clusters generically degenerate into five.
Physically, there are still seven clusters of distinguishable modes (eigenvectors), it is just that in the absence of viscosity the two vortical clusters are indistinguishable from \(\lambda\approx 0\) in an eigenvalue plot.
Importantly, \cref{fig:compEigs_Sp_cD0cV0} illustrates that the patch scheme preserves \(\Re(\lambda)=0\) in the case of ideal waves.

\subsection{Consistent accuracy of polynomial coupling}

\begin{figure}
  \caption{\label{fig:compEigs_Sq2468}%
     Macroscale eigenvalues~\(\lambda_p\) (\(\Re(\lambda_p) > -0.001\)) of polynomial patch scheme (\(N = 10\), \(n = 6\), \(r=0.3\)) agree better with eigenvalues~\(\lambda_\mu\) of the full-domain model as the order~\(p\) of the polynomial coupling increases.
    }
  \centering
  \begin{subfigure}[b]{0.489\textwidth}
    \centering
    \caption{\label{fig:compEigs_Sq2}%
      square-p2 patch scheme
      }
    \includegraphics[scale=0.95]{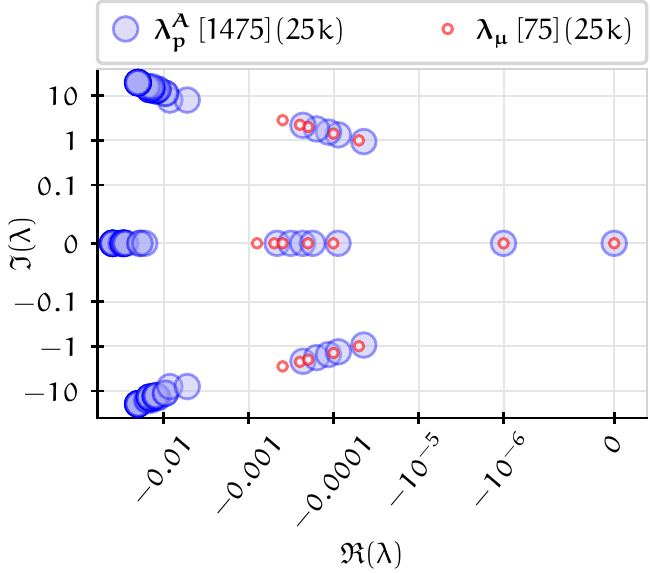}
  \end{subfigure}
  \begin{subfigure}[b]{0.489\textwidth}
    \centering
    \caption{\label{fig:compEigs_Sq4}%
      square-p4 patch scheme
      }
    \includegraphics[scale=0.95]{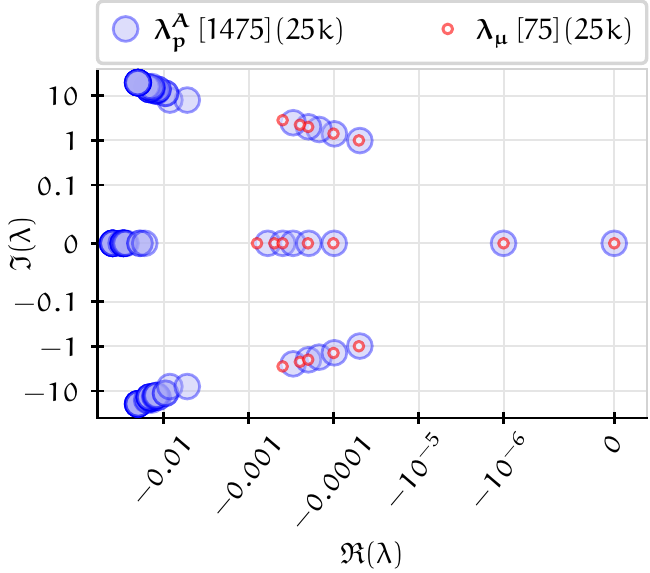}
  \end{subfigure}
  \begin{subfigure}[b]{0.489\textwidth}
    \centering
    \caption{\label{fig:compEigs_Sq6}%
      square-p6 patch scheme
      }
    \includegraphics[scale=0.95]{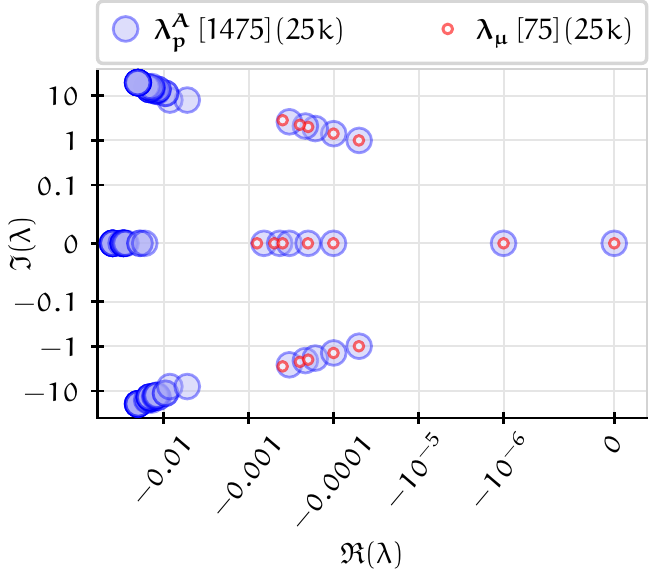}
  \end{subfigure}
  \begin{subfigure}[b]{0.489\textwidth}
    \centering
    \caption{\label{fig:compEigs_Sq8}%
      square-p8 patch scheme
      }
    \includegraphics[scale=0.95]{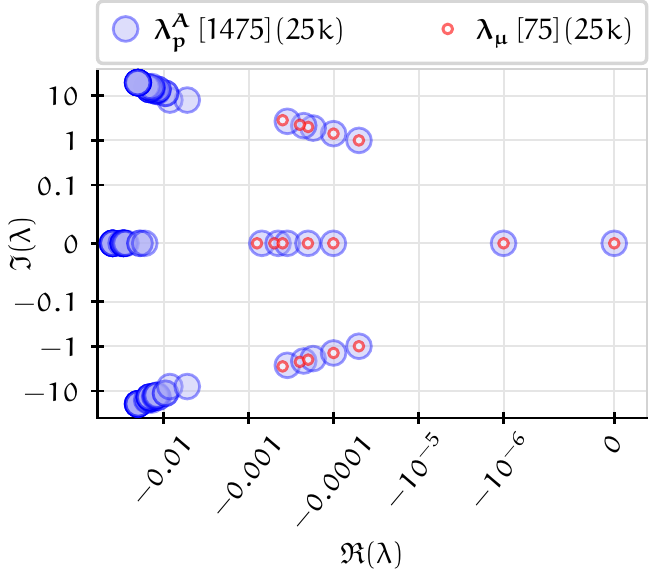}
  \end{subfigure}
\end{figure}

\Cref{fig:compEigs_Sq2468} visually compares the eigenvalues~\(\lambda_p\) of the polynomial patch schemes (\cref{ssc:plynmlCplng}, (\(p=2,4,6,8\)) with the eigenvalues~\(\lambda_\mu\) of the full-domain model~\eqref{eqs:FDEs_mN_gLnrWve} for weakly damped linear waves with~\(c_D = 10^{-6}\) and \(c_V = 10^{-4}\). 
\Cref{fig:compEigs_Sq2468} shows that with increasing order~\(p\), the macroscale eigenvalues~\(\lambda_{p,M}\) (\(\Re(\lambda_{p}) > -0.001\)) of the polynomial schemes improve agreement with the corresponding eigenvalues~\(\lambda_\mu\) of the full-domain model.
That is, \emph{increasing interpolation order~\(p\) of the patch coupling improves accuracy}.
\Cref{ssc:plynmlPtchSchmsAreCnsstnt_gLnrWve} confirms this trend quantitatively for a wide range of grid and physical parameters.

\begin{figure}
  \caption{\label{fig:compEigs_Sq4_Nn}%
    {Increasing macro-grid intervals~\(N\) (left to right) improves accuracy} of polynomial patch schemes, here for square-p4. 
    Increasing~\(n\) (top to bottom) does not affect accuracy, but does increase the number of microscale modes (\(\Re(\lambda_p) < -0.01\)). 
    Patch ratio~\(r=0.1\).
    }
  \centering
\begin{tabular}{@{}cc@{}}
  \begin{subfigure}[b]{0.489\textwidth}
    \centering
    \caption{\label{fig:compEigs_Sq4_N6n6}%
      \(N = 6\), \(n = 6\), \(\delta = \pi/90\). % \delta = pi /(N*n/(4*r))
      }
    \includegraphics[scale=0.95]{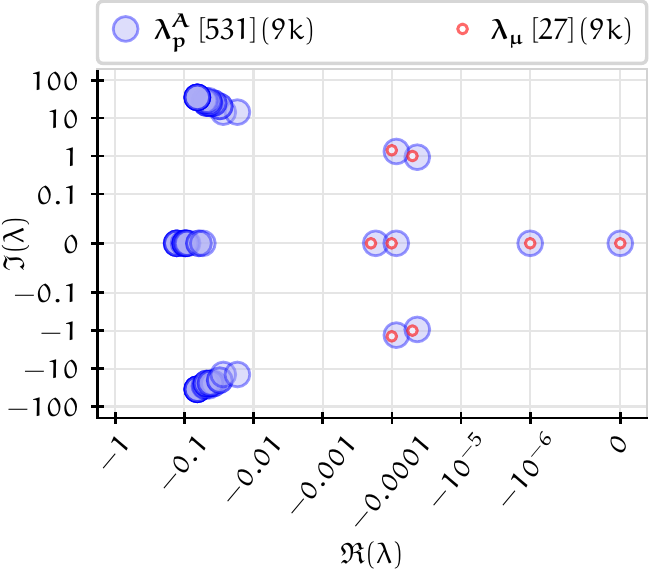}
  \end{subfigure}
  &
  \begin{subfigure}[b]{0.489\textwidth}
    \centering
    \caption{\label{fig:compEigs_Sq4_N10n6}%
      \(N = 10\), \(n = 6\), \(\delta = \pi/150\).
      }
    \includegraphics[scale=0.95]{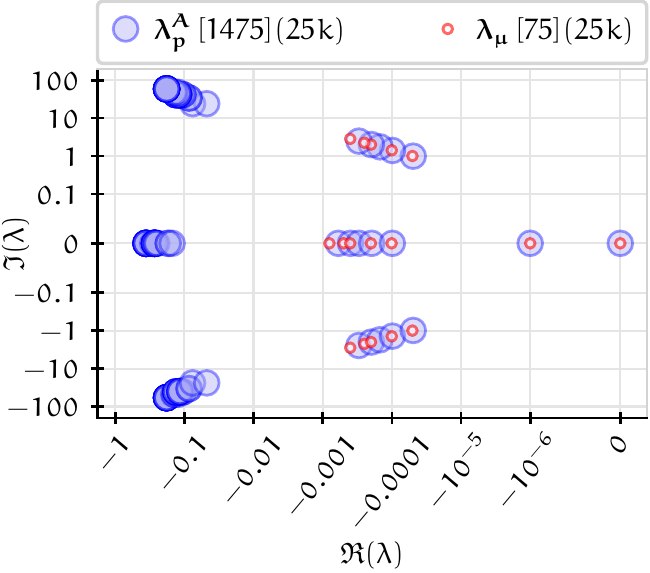}
  \end{subfigure}
  \\
  \begin{subfigure}[b]{0.489\textwidth}
    \centering
    \caption{\label{fig:compEigs_Sq4_N6n10}%
      \(N = 6\), \(n = 10\), \(\delta = \pi/150\).
      }
    \includegraphics[scale=0.95]{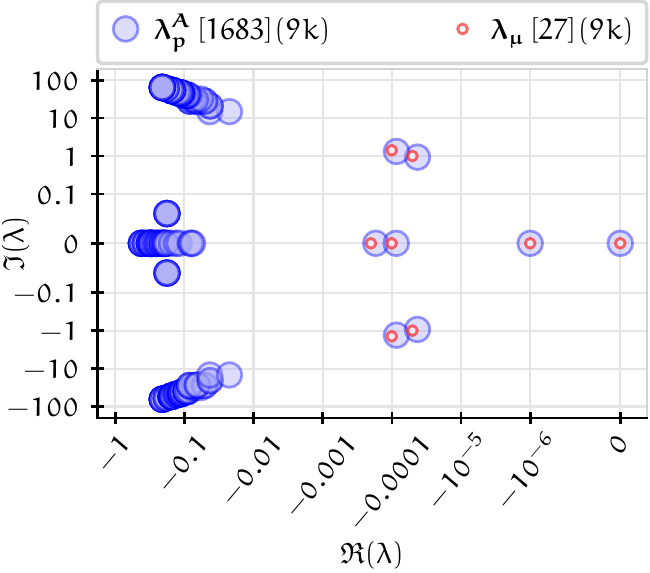}
  \end{subfigure}
  &
  \begin{subfigure}[b]{0.489\textwidth}
    \centering
    \caption{\label{fig:compEigs_Sq4_N10n10}%
      \(N = 10\), \(n = 10\), \(\delta = \pi/250\).
      }
    \includegraphics[scale=0.95]{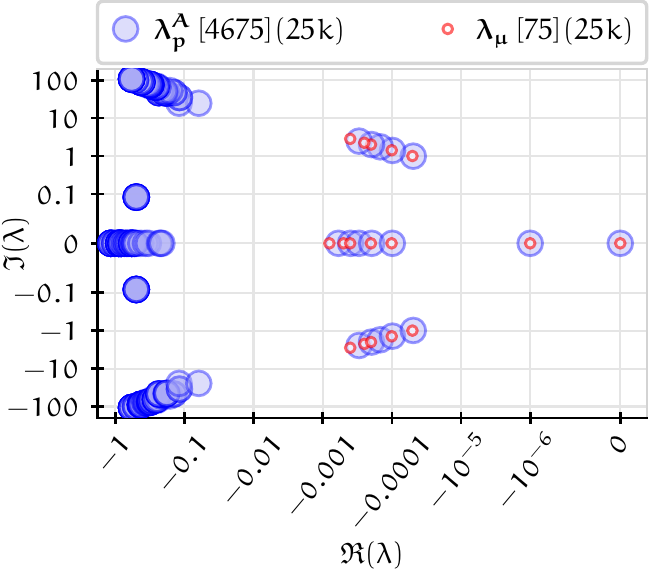}
  \end{subfigure}
\end{tabular}
\end{figure}

Increasing the number of macro-grid intervals~\(N\) increases the accuracy of the polynomial patch schemes.
For all four combinations of the grid parameters~\(N, n \in \{6, 10\}\), all with patch ratio~\(r = 0.1\), \Cref{fig:compEigs_Sq4_Nn} visually compares the eigenvalues~\(\lambda_p\) of the square-p4 polynomial patch scheme with the eigenvalues~\(\lambda_\mu\) of the full-domain weakly damped linear waves~\cref{eqs:FDEs_mN_gLnrWve} with~\(c_D = 10^{-6}\) and \(c_V = 10^{-4}\).
\Cref{fig:compEigs_Sq4_Nn} indicates that upon increasing~\(N\), the macroscale eigenvalues~\(\lambda_{p,M}\) of the polynomial schemes agree better with the corresponding eigenvalues~\(\lambda_\mu\) of the full-domain eigenvalues.
On the other hand, increasing the number of sub-patch micro-grid intervals~\(n\) does not affect the patch scheme's macroscale accuracy, but increases computation cost by increasing the number of microscale modes (\(\Re(\lambda_p) < -0.001\) in \cref{fig:compEigs_Sq2468}).
Similarly, comparing \cref{fig:compEigs_Sq4} for \(r = 0.3\) and \cref{fig:compEigs_Sq4_N10n6} for \(r = 0.1\) indicates that changing the patch ratio~\(r\) (keeping~\(N,n\) constant) also does not have any effect on the patch scheme accuracy.

\section{The schemes are not sensitive to numerical roundoff errors}
% ++++++++++++++++++++++++++++++++++++++++++++++++++++++++++++++++++++
\label{sec:ptchSchmsAreNtSnstveToNmrclErrs_gLnrWve}

From both qualitative and quantitative arguments, this section shows that, for weakly damped linear waves, the staggered patch schemes are not sensitive to numerical roundoff errors, and hence are suitable for practical use
\parencite[e.g.,][]{Goldberg1991_WhtEvryScntstShldKnwAbtFltngPntArthmtc}.
This insensitivity to roundoff errors empowers the scheme to accurately resolve complex physics over a wide range of length scales.
The insensitivity also allows using advanced optimisations in finite precision computing, such as \textsc{simd} loop reordering, fast math mode, and treating subnormal numbers as zeros \parencite{Goldberg1991_WhtEvryScntstShldKnwAbtFltngPntArthmtc}.

Recall from \cref{sec:ptchSchmsAreAcrte_gLnrWve} that here errors are best seen in the eigenvalues of the system.
Hence to explore roundoff errors, this section compares the eigenvalues~\(\lambda_p^A\) of the analytically derived Jacobian \parencite[\S4.1]{Divahar2022_StgrdGrdsFrMltdmnsnlMltscleMdlng} with the eigenvalues~\(\lambda_p^N\) of the numerically computed Jacobian.
We define the \emph{microscale and macroscale roundoff errors} for the staggered patch scheme eigenvalues as
\begin{equation} \label{eqs:errNum_defs}
    \varepsilon_\mu := \max_i |\lambda_{p\,\mu,i}^N - \lambda_{p\,\mu,i}^A|
      \,,\qquad
    \varepsilon_M := \max_i |\lambda_{p\,M,i}^N - \lambda_{p\,M,i}^A|\,,
\end{equation}
for eigenvalues ordered in index~\(i\) to correspond appropriately.
Following the method of \textcite[\S3.3]{Divahar2022_AcrteMltiscleSmltnOfWveLkeSystms}, we group the patch scheme numerical eigenvalues~\(\lambda_p^N\) wavenumber-wise, establish the pair-wise correspondence between the analytical and numerical eigenvalues~(\(\lambda_p^A,\lambda_p^N\)), and separate the eigenvalues~\(\lambda_p^N\) into microscale and macroscale eigenvalues~(\(\lambda_{p,\mu}^N,\lambda_{p,M}^N\)) respectively.
When the errors~\(\varepsilon_\mu\) and~\(\varepsilon_M\) are negligibly small, then the patch scheme is not sensitive to roundoff errors.

\begin{table}
  \centering
  \caption{\label{tbl:prmtrs_errNum_gLnrWve_alSchms}%
    \Cref{sec:ptchSchmsAreNtSnstveToNmrclErrs_gLnrWve} studies the sensitivity of patch scheme eigenvalues to roundoff errors for all the~\(1944\) combinations of these parameters. % 9*3*2*4 + 4*9*6*2*4 = 1944
  }
  \begin{tabular}{ll}
  \hline
  Patch schemes &
    spectral, square-p2, square-p4, square-p6,  square-p8
    \\
  Dissipation &
    \(c_D \in \{0, 10^{-6}, 0.001\}\),\quad
    \(c_V \in \{0, 10^{-4}, 0.01\}\)
    \\
  Macro-grid &
    \(N \in \{6, 10, 14\}\) for spectral interpolation,\\
    &
    \(N \in \{6, 10, 14, 18, 22, 26\}\) for polynomial interpolation
    \\
  Micro-grid &
    \(n \in \{6, 10\}\)
    \\
  Patch ratio & 
    \(r \in \{0.0001, 0.001, 0.01, 0.1\}\)
  \\\hline
  \end{tabular}
\end{table}
We computed the sensitivity of the patch scheme eigenvalues to roundoff errors for all the~\(1944\) combinations of the listed parameters in \cref{tbl:prmtrs_errNum_gLnrWve_alSchms}. % 9*3*2*4 + 4*9*6*2*4 = 1944
\begin{figure}
  \centering
  \caption{\label{fig:eNum_vs_Nr_micro_macro_gLnrWve}%
    For all five patch schemes, variation of peak microscale and macroscale roundoff errors with the number of macro- and micro-grid intervals~\(N\), \(n\) and patch ratio~\(r\).
  }
  \includegraphics[scale=0.88]{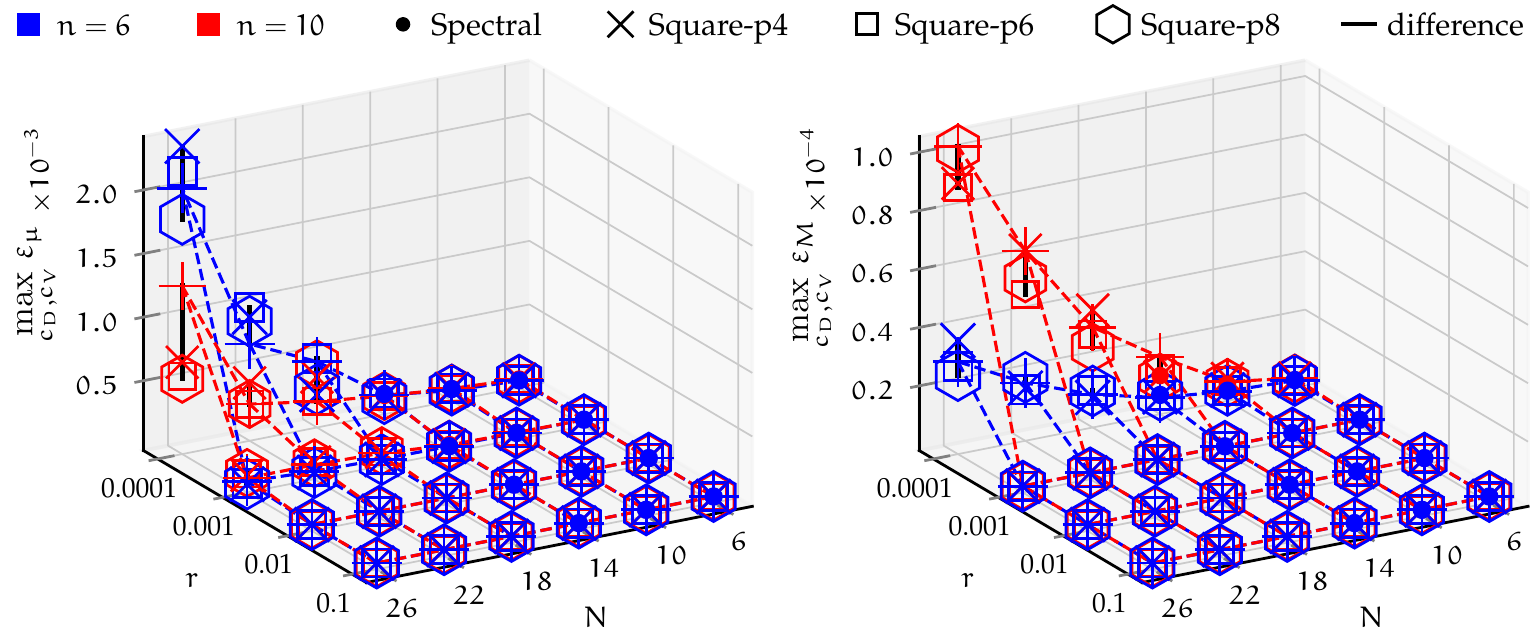}
\end{figure}%
For an overall summary of a patch scheme on a patch grid with~\(N\), \(n\), and~\(r\), we define the \emph{peak microscale and macroscale roundoff errors} as~\(\max_{c_D,c_V} \varepsilon_\mu\) and~\(\max_{c_D,c_V} \varepsilon_M\) respectively, over all the nine combinations of the coefficients~\(c_D,c_V\) in \cref{tbl:prmtrs_errNum_gLnrWve_alSchms}.
\Cref{fig:eNum_vs_Nr_micro_macro_gLnrWve,fig:eNum_vs_Nr_micro_macro_gLnrWve} shows that \emph{nonnegligible roundoff errors arise only for very small sub-patch micro-grid spacing~\(\delta \lesssim 10^{-5}\) (i.e., small~\(r\) and large~\(N,n\)).}
The following are some key observations from \cref{fig:eNum_vs_Nr_micro_macro_gLnrWve,fig:eNum_vs_Nr_micro_macro_gLnrWve} (\textcite{Divahar2022_AcrteMltiscleSmltnOfWveLkeSystms} discusses many more details).
\begin{itemize}
  \item \label{itm:macroNumErrSmall} For each set of patch grid parameters~\(N\), \(n\), and~\(r\), the peak macroscale roundoff errors~\(\max_{c_D,c_V} \varepsilon_M\) are about ten to thousand times smaller than the peak microscale roundoff errors~\(\max_{c_D,c_V} \varepsilon_\mu\).
  That is, \emph{the macroscale eigenvalues which are of primary interest are less sensitive to roundoff errors than the microscale eigenvalues}.
  \item Both the microscale and macroscale peak roundoff errors~\(\max_{c_D,c_V} \varepsilon_\mu\) and~\(\max_{c_D,c_V} \varepsilon_M\), monotonically increase with increasing number of macro-grid intervals~\(N\) and decreasing patch ratio~\(r\).
  The roundoff errors also increase with increasing number of sub-patch micro-grid intervals~\(n\) (blue and red in \cref{fig:eNum_vs_Nr_micro_macro_gLnrWve,fig:eNum_vs_Nr_micro_macro_gLnrWve}), except the off trend for~\(\max_{c_D,c_V} \varepsilon_\mu\) in \cref{fig:eNum_vs_Nr_micro_macro_gLnrWve} for \(N \in \{22,26\}\), \(r = 0.0001\).
  For a staggered patch grid, increasing~\(N\), decreasing~\(r\), and increasing~\(n\), all these lead to decreasing sub-patch micro-grid spacing~\(\delta = 2 L r / (N n)\).
  For example, for \(r = 0.001\), \(N=26\), \(n=10\), sub-patch micro-grid spacing~\(\delta \approx 5 \cdot 10^{-5}\)
  Thus, \emph{nonnegligible roundoff errors arise only for very small sub-patch micro-grid spacing~\(\delta \lesssim 10^{-5}\) (i.e., small~\(r\) and large~\(N,n\)).}
  \item In general, except~\(\varepsilon_\mu\) for \(N \gtrsim 22\) and \(r=0.0001\), the roundoff errors of all five patch schemes are roughly the same.
  That is, the roundoff errors do not have a strong dependence on the specific patch scheme.
  If the roundoff errors were due to the patch scheme, then the roundoff errors must also depend on the specific patch scheme, showing a clear trend.
  The lack of such trends, among the patch schemes with different amounts of numerical computations, indicates that \emph{the roundoff errors are not due to the patch schemes}.
\end{itemize}

Thus, except for very small sub-patch micro-grid spacing~\(\delta\), the roundoff errors~\(\varepsilon_\mu\) and~\(\varepsilon_M\) are small, and the roundoff errors do not depend on the specific patch scheme.
Hence, the staggered patch schemes are effectively not sensitive to roundoff errors.

\section{Staggered patch schemes are stable}
% ++++++++++++++++++++++++++++++++++++++++++++++++++++++++++++++++++++
\label{sec:ptchSchmsAreStble_gLnrWve}

This section demonstrates the stability of the \emph{spectral and square-p staggered patch schemes}, for a range of patch scheme parameters, for nine combinations of the physical parameters~\(c_D, c_V\).

\begin{table}
  \centering
  \caption{%
    \Cref{sec:ptchSchmsAreStble_gLnrWve} studies the stability of the staggered patch schemes using the eigenvalues for all the~\(4\,374\) combinations of these parameters.
  }
  \label{tbl:prmtrs_stblty_gLnrWve_alSchms}
  \renewcommand{\arraystretch}{1.1}
  \begin{tabular}{lp{0.65\textwidth}}
  \hline
  Patch schemes &
    spectral, square-p2, square-p4, square-p6, and square-p8
    \\
  Drag coefficient &
    \(c_D \in \{0, 10^{-6}, 0.001\}\)
    \\
  Viscous coefficient &
    \(c_V \in \{0, 10^{-4}, 0.01\}\)
    \\
  Macro-grid &
    \(N \in \{6, 10, 14\}\) for spectral scheme,
    \\ &
    \(N \in \{6, 10, 14, 18, 22, 26\}\) for polynomial schemes.
    \\
  Micro-grid &
    \(n \in \{6, 10\}\)
    \\
  Patch ratio & 
    \(r \in \{0.0001, 0.001, 0.01, 0.1\}\)
  \\\hline
  \end{tabular}
\end{table}%
A patch scheme may potentially be unstable due to either the macroscale modes or the microscale modes.
So we computed the maximum real parts of the numerical eigenvalues of the five staggered patch schemes, separately for the microscale and macroscale modes (\(\max \Re(\lambda_{p\,\mu}^N)\) and \(\max \Re(\lambda_{p\,M}^N)\)), for the~\(4\,374\) combinations of the parameters in \cref{tbl:prmtrs_stblty_gLnrWve_alSchms}.
\begin{table}
  \centering
  \caption{%
    Overall maximum real parts of the microscale and macroscale eigenvalues (\(\lambda_{p\,\mu}^N\), \(\lambda_{p\,M}^N\)) over the five patch schemes and all the combinations of \(c_D\), \(c_V\) and \(N\) in \cref{tbl:prmtrs_stblty_gLnrWve_alSchms}, for different number of sub-patch intervals~\(n\) and patch ratios~\(r\).
  }
  \label{tbl:peakMaxRe_gLnrWve}
\begin{math}
  \renewcommand{\arraystretch}{1.1}
  \begin{array}{llrrrr}
  {} & {} & \multicolumn{4}{c}{\text{Patch ratio}~r} \\
  {} &
    {} &
    0.0001 &
    0.001 &
    0.01 &
    0.1
    \\ \cline{3-6}
  \multirow{2}{*}{Overall  \(\max \Re(\lambda_{p\,\mu}^N)\)} &
    n=6 &
    2\cdot 10^{-6} &
    2\cdot 10^{-8} &
    2\cdot 10^{-10} &
    2\cdot 10^{-12}
    \\
  {} &
    n=10 &
    {\color{black}5\cdot 10^{-6}} &
    {\color{black}7\cdot 10^{-8}} &
    {\color{black}6\cdot 10^{-10}} &
    {\color{black}7\cdot 10^{-12}}
    \\\hline
  \multirow{2}{*}{Overall  \(\max \Re(\lambda_{p\,M}^N)\)} &
    n=6 &
    7\cdot 10^{-6} &
    8\cdot 10^{-9} &
    3\cdot 10^{-11} &
    10^{-12}
    \\
  {} &
    n=10 &
    {\color{black}3\cdot 10^{-5}} &
    {\color{black}2\cdot 10^{-8}} &
    {\color{black}2\cdot 10^{-10}} &
    {\color{black}3\cdot 10^{-12}}
  \\ \cline{3-6}
  \end{array}
\end{math}
\end{table}%
\Cref{tbl:peakMaxRe_gLnrWve} lists the overall maximum real parts. 
It shows that for moderately small patch ratios \(r \in \{0.01, 0.1\}\), both the microscale and macroscale eigenvalues of the patch schemes have the maximum real parts less than about~\(6 \cdot 10^{-10}\).
Thus, \emph{for moderately small patch ratios~\(r \gtrsim 0.01\), all five patch schemes are stable}.

For smaller patch ratios~\(r \in \{0.0001, 0.001\}\), \Cref{tbl:peakMaxRe_gLnrWve} shows that some of the patch schemes have  maximum real parts up to about~\(10^{-5}\) which correspond to only a few specific combinations of the parameters of the patch grid~(\(N,n,r\)) and the physical system~(\(C_D,C_V\)).
All such combinations of grid parameters correspond to a very small micro-grid spacing of \(\delta \lesssim 10^{-5}\) for which the microscale computations incur increased roundoff error.  
\textcite[\S3.4]{Divahar2022_AcrteMltiscleSmltnOfWveLkeSystms} shows with detailed evidence that such nonnegligible real parts are due to eigenvalue computation being affected by roundoff errors, either due to the many repeated near-zero microscale eigenvalues (as in \cref{fig:compEigs_Sp_cD0cV0,fig:compEigs_Sp_cD1e-6cV0}), or due to the inherent sensitivity of the microscale model affecting the accurate computation of the three near-zero macroscale eigenvalues (as in \cref{fig:compEigs_Sp_cD0cV1e-4,fig:compEigs_Sp_cD1e-6cV1e-4}).
Thus such nonnegligible real parts of about~\(10^{-5}\) are not due to the patch scheme: \emph{the patch schemes developed herein are stable}.

\section{The schemes are consistent with the given microscale model}
% ++++++++++++++++++++++++++++++++++++++++++++++++++++++++++++++++++++
\label{sec:ptchSchmsAreCnsstnt_gLnrWve}

This section shows that the staggered patch schemes are \emph{consistent} with the given microscale model.
A computational model is usually called consistent when the discretized equations, such as~\cref{eqs:FDEs_mN_gLnrWve}, approach to the corresponding \textsc{pde}s~\cref{eqs:PDEs_gLnrWve}, as the micro-grid spacing~\(\delta \to 0\) \cite[e.g.,][p.~34]{Ferziger2020_CmpttnlMthdsFrFldDynmcs}.
But the goal of our multiscale staggered patch scheme~\cref{eqn:dynSys_pN_gnrcWve} is to accurately represent the \emph{macroscale} waves of the corresponding \emph{discrete} full-domain microscale model~\cref{eqs:FDEs_mN_gLnrWve}.
Hence we \emph{define a patch scheme to be consistent} when the macroscale characteristics of the patch scheme \cref{eqn:dynSys_pN_gnrcWve} approach to the corresponding macroscale characteristics of the full-domain microscale model~\cref{eqs:FDEs_mN_gLnrWve} with decreasing patch spacing~\(\Delta\).

This section establishes the consistency of the staggered patch schemes by showing that the macroscale eigenvalues~\(\lambda^N_{p\,M}\) of the patch schemes converge to the macroscale eigenvalues of the corresponding full-domain microscale model with decreasing patch spacing~\(\Delta\).
The eigenvalue spectra in \cref{sec:ptchSchmsAreAcrte_gLnrWve} show that the staggered patch scheme macroscale eigenvalues~\(\lambda_{p\,M}^A\) (e.g., clusters~\(1\)--\(4\) in \cref{fig:compEigs_Sp_cD1e-6cV1e-4}) have  similar qualitative structure, and are visually close, to the corresponding macroscale eigenvalues~\(\lambda_\mu\) of the fine-grid full-domain microscale model (by varying degrees depending upon the specific staggered patch scheme,~\(N\), \(n\), and~\(r\)).
To numerically quantify the discrepancy between the macroscale eigenvalues~\(\lambda_{p\,M}^N\) and~\(\lambda_\mu\), we define the \emph{eigenvalue error} for the \emph{macroscale} wavenumber~\((k_x,k_y)\) as
\begin{equation} \label{eqn:eigErr_gLnrWve}
  \epsilon^{k_x,k_y} := \| \underline{\lambda}^N_{p\,M}(k_x,k_y) - \underline{\lambda}_\mu(k_x,k_y) \| \big/ \| \underline{\lambda}_\mu(k_x,k_y) \|,
\end{equation}
where~\(\|\argdot\|\) is the Euclidean norm of the three element complex vectors of eigenvalues~\(\underline{\lambda}^N_{p\,M}\) and \( \underline{\lambda}_\mu\) (three macroscale eigenvalues for each macroscale wavenumber).

\begin{table}
  \centering
  \caption{%
    \Cref{sec:ptchSchmsAreCnsstnt_gLnrWve} establishes the consistency of the patch schemes using eigenvalues for all the~\(2\,160\) combinations of these parameters.
    The inter-patch spacing, here \(\Delta=2\pi/N\), decreases with increasing~\(N\).  
  }
  \label{tbl:prmtrs_cnvrgnce_gLnrWve}
  \renewcommand{\arraystretch}{1.1}
  \begin{tabular}{lp{0.65\textwidth}}
  \hline
  Patch schemes &
    spectral, square-p2, square-p4, square-p6, and square-p8
    \\
  Drag coefficient &
    \(c_D \in \{0, 10^{-6}, 0.001\}\)
    \\
  Viscous coefficient &
    \(c_V \in \{0, 10^{-4}, 0.01\}\)
    \\
  Macro-grid &
    \(N \in \{6, 10, 14, 18, 22, 26\}\)
    \\
  Micro-grid &
    \(n \in \{6, 10\}\)
    \\
  Patch ratio & 
    \(r \in \{0.0001, 0.001, 0.01, 0.1\}\)
  \\[0.5ex]\hline
  \end{tabular}
\end{table}

To assess the patch scheme consistency in this section (i.e., eigenvalue convergence), we compute the three eigenvalue errors~\(\epsilon^{1,0}\), \(\epsilon^{1,1}\) and~\(\epsilon^{2,1}\), for the patch schemes corresponding to the three macroscale (angular) wavenumbers \((k_x,k_y) \in \{(1,0), (1,1), (2,1)\}\) over the \(2\pi \times 2\pi\) non-dimensional domain.
For example, the smallest wavenumber~\((1,0)\) corresponds to the largest wavelength of~\((2\pi, 0)\) over the chosen \(2\pi \times 2\pi\) domain.
We computed these errors for all~\(2\,160\) combinations of the parameters listed in \cref{tbl:prmtrs_cnvrgnce_gLnrWve}.

Computing the three element vector of eigenvalues~\(\underline{\lambda}_\mu(k_x,k_y)\) in the eigenvalue error~\cref{eqn:eigErr_gLnrWve} is straightforward using~\cref{eqn:eig_muA}.
On the other hand, among the numerical eigenvalues~\(\lambda_{p}^N\), finding which three eigenvalues correspond to the three eigenvalues~\(\underline{\lambda}_\mu(k_x,k_y)\) (for the same macroscale wavenumber), is not straightforward.
A heuristic method (based only on the eigenvalues without using the eigenvectors) by \textcite[\S3.3]{Divahar2022_AcrteMltiscleSmltnOfWveLkeSystms} classifies the patch scheme eigenvalues~\(\lambda_{p}^N\) wavenumber-wise, and separates the eigenvalues into microscale and macroscale patch scheme eigenvalues~\(\lambda_{\mu\,M}^N, \lambda_{p\,M}^N\).
The same method also establishes the correspondence between the three-element vector of full-domain microscale eigenvalues~\( \underline{\lambda}_\mu\) and the patch scheme eigenvalue vector~\(\underline{\lambda}_{p\,M}^N\) in the definition~\cref{eqn:eigErr_gLnrWve}.

\subsection{Spectral patch scheme is uniformly consistent}
% --------------------------------------------------------------------
\label{ssc:spctrlPtchSchmsAreCnsstnt_gLnrWve}

With the highly accurate global spectral interpolation (\cref{ssc:spctrlCplng}), the \emph{spectral patch scheme resolves the macroscale modes exactly} to within roundoff errors, irrespective of the number of patch spacing~\(\Delta\) (e.g., the complex plane eigenvalue plot \cref{fig:compEigs_Sp_cDcV} of \cref{ssc:SpctrlCplngIsHghlyAcrte}).
That is, the spectral patch scheme is uniformly consistent with the given microscale model without any dependence on the patch spacing~\(\Delta\).

\begin{table}
  \centering
  \caption{%
    Logarithm of the maximum eigenvalue error~\( \log_{10} \max_N \epsilon^{k_x,k_y}\) for the spectral staggered patch scheme over the six different number of macro-grid intervals~\(N\) in \cref{tbl:prmtrs_cnvrgnce_gLnrWve}, and over wavenumbers \((k_x,k_y)\in\{(1,0),(1,1),\,(2,1)\}\).
    Red colour highlights \(\epsilon^{k_x,k_y} > 2\cdot10^{-8}\). 
  }
  \label{tbl:eigErr_gLnrWve_cDcV_Spctrl}
  \renewcommand{\arraystretch}{1.1}
\begin{math}
\begin{array}{lrrrrrrrr}
  {} & \multicolumn{8}{c}{\text{Patch ratio }r} \\
  {} &
    \multicolumn{2}{c}{0.0001}&
    \multicolumn{2}{c}{0.001}&
    \multicolumn{2}{c}{0.01}&
    \multicolumn{2}{c}{0.1}
    \\ \cline{2-9}
  c_D, c_V \ \backslash\ n &
    6 &
    10 &
    6 &
    10 &
    6 &
    10 &
    6 &
    10
    \\ \hline 
  0, 0 &
      -10 &
     -9.3 &
      -12 &
      -11 &
      -12 &
      -12 &
      -12 &
      -12
    \\
  0, 0.0001 &
     -9.4 &
     -8.2 &
      -11 &
      -10 &
      -12 &
      -12 &
      -12 &
      -12
    \\
  0, 0.01 &
    {\color{Red} -5.6} &
    {\color{Red} -5.3} &
     -8.5 &
       -8 &
      -11 &
      -10 &
      -12 &
      -12
    \\ \hline 
  10^{-6}, 0 &
     -9.9 &
     -9.1 &
      -11 &
      -10 &
      -12 &
      -12 &
      -12 &
      -12
    \\
  10^{-6}, 0.0001 &
       -9 &
     -8.3 &
      -11 &
      -10 &
      -12 &
      -12 &
      -12 &
      -12
    \\
  10^{-6}, 0.01 &
    {\color{Red} -5.5} &
    {\color{Red} -5.1} &
     -8.6 &
     -7.9 &
      -11 &
      -10 &
      -12 &
      -12
    \\ \hline
  0.001, 0 &
      -10 &
     -9.2 &
      -11 &
      -11 &
      -12 &
      -11 &
      -12 &
      -12
    \\
  0.001, 0.0001 &
       -9 &
     -8.2 &
      -11 &
      -10 &
      -12 &
      -12 &
      -12 &
      -12
    \\
  0.001, 0.01 &
    {\color{Red} -5.7} &
    {\color{Red} -4.9} &
     -8.4 &
       -8 &
      -11 &
      -10 &
      -12 &
      -12
  \\
  \hline
  \end{array}
\end{math}  
\end{table}

\begin{figure}  
\caption{\label{fig:iD2_iV2_eVsD_10_pN_m1dNE_Spectral_cD0p001_cV0p01_n10}%
    Worst-case eigenvalue errors~\(\epsilon^{1,0}\) of the spectral staggered patch scheme, for \(c_D=0.001\), \(c_V = 0.01\), \(n=10\), for different patch spacings~\(\Delta\) and patch ratio~\(r\).  Plots for~\(\epsilon^{1,1}\) and~\(\epsilon^{2,1}\) are equivalent.}
  
  \centering
  \includegraphics{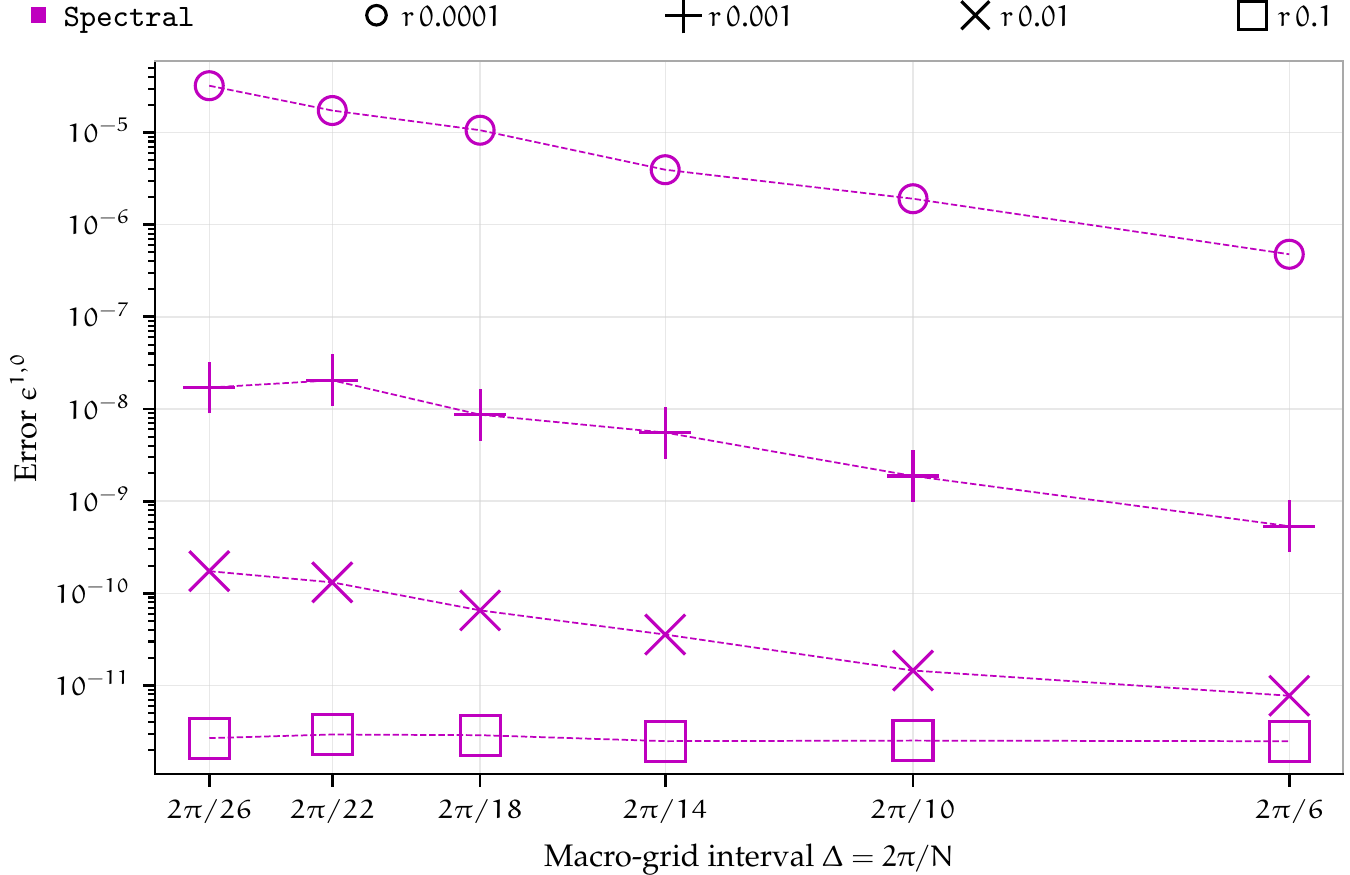}
\end{figure}

\Cref{tbl:eigErr_gLnrWve_cDcV_Spctrl} shows the maximum eigenvalue errors of~\(\epsilon^{1,0}\), \(\epsilon^{1,1}\) and~\(\epsilon^{2,1}\), over the six different number of macro-grid intervals~\(N\) in \cref{tbl:prmtrs_cnvrgnce_gLnrWve}.
For the worst case (i.e., largest eigenvalue errors) among the combinations of the parameters in in \cref{tbl:prmtrs_cnvrgnce_gLnrWve}, \cref{fig:iD2_iV2_eVsD_10_pN_m1dNE_Spectral_cD0p001_cV0p01_n10} plots the variation of the eigenvalue errors~\(\epsilon^{1,0}\) with the patch spacing \(\Delta = 2\pi/N\).
\Cref{tbl:eigErr_gLnrWve_cDcV_Spctrl,fig:iD2_iV2_eVsD_10_pN_m1dNE_Spectral_cD0p001_cV0p01_n10} together show that, except for the combination of the small patch ratio~\(r \lesssim 0.01\), small patch spacing~\(\Delta \lesssim 2\pi/18\) and large viscosity~\(c_V = 0.01\), all the three eigenvalue errors are about~\(10^{-8}\) or smaller.
\textcite[\S3.6.1]{Divahar2022_AcrteMltiscleSmltnOfWveLkeSystms} reports more details of these small errors.
\cref{sec:ptchSchmsAreNtSnstveToNmrclErrs_gLnrWve} indicates that the cases of larger eigenvalue errors (larger than~\(10^{-8}\)) are due to roundoff errors.
This small error shows that the spectral patch scheme is uniformly consistent without any dependence on the patch spacing.

\subsection{Polynomial patch schemes are consistent}
% --------------------------------------------------------------------
\label{ssc:plynmlPtchSchmsAreCnsstnt_gLnrWve}

For the patch coupling by polynomial interpolation, comparing the subplots of \cref{fig:compEigs_Sq4_Nn} on left (\(N=6\)) with the subplots on right (\(N=10\)) indicates that decreasing the patch spacing~\(\Delta\) (increasing~\(N\)) improves the patch scheme accuracy.
We confirmed and characterised this increasing accuracy with decreasing~\(\Delta\) (i.e., consistency) by exploring the macroscale eigenvalue errors~\(\epsilon^{1,0}\), \(\epsilon^{1,1}\) and~\(\epsilon^{2,1}\) (defined by~\cref{eqn:eigErr_gLnrWve}) for all the~\(1\,728\) combinations of the parameters in \cref{tbl:prmtrs_cnvrgnce_gLnrWve} for the polynomial patch schemes.

\begin{figure}
\caption{\label{fig:iD2_iV2_eVsD_10_pN_m1dNE_Square_cD0p001_cV0p01_n10}
    Worst-case convergence of macroscale eigenvalues with patch spacing~\(\Delta\) (\((k_x,k_y)=(1,0)\), \(c_D=0.001\), \(c_V = 0.01\), and \(n=10\)) for the four polynomial patch schemes (interpolation orders~\(p \in \{2, 4, 6, 8\}\)) and various patch ratios~\(r\). Solid lines are the power law  fit \({\epsilon}^{1,0} = 0.3333 \cdot (0.7 \cdot \Delta)^p \).   }
  \centering
  \includegraphics{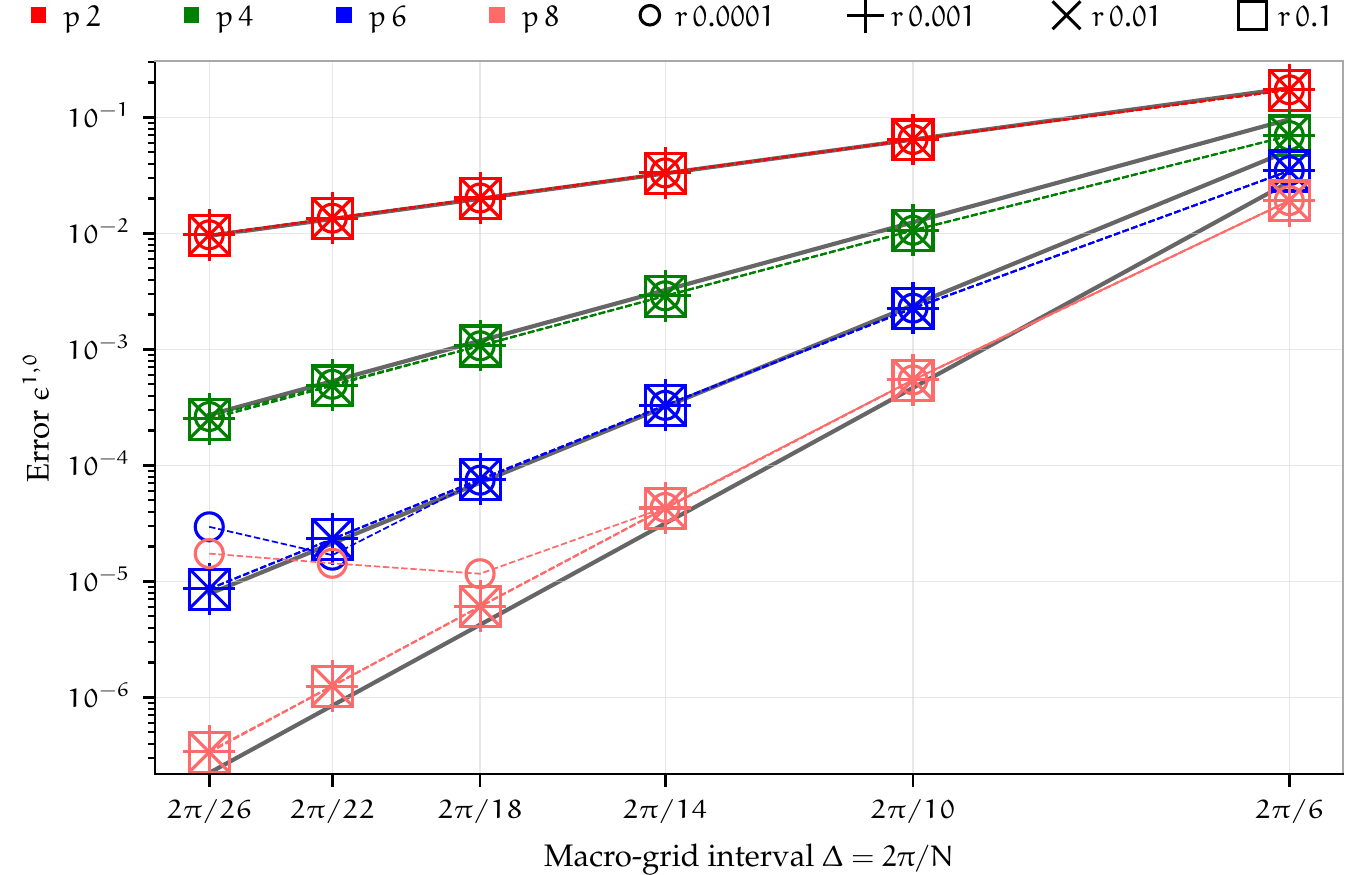}
\end{figure}

\cref{fig:iD2_iV2_eVsD_10_pN_m1dNE_Square_cD0p001_cV0p01_n10} shows an example case of the variation of macroscale eigenvalue error with patch spacing~\(\Delta\).
As the patch spacing~\(\Delta\) decreases the plotted error decreases and hence shows that \emph{the polynomial patch schemes are consistent with the underlying given microscale model}.
The case plotted in \cref{fig:iD2_iV2_eVsD_10_pN_m1dNE_Square_cD0p001_cV0p01_n10} is the `worst-case' example over all the parameters explored (\cref{tbl:prmtrs_cnvrgnce_gLnrWve}).
\textcite[\S3.6.2]{Divahar2022_AcrteMltiscleSmltnOfWveLkeSystms} reports more cases and details of the macroscale eigenvalue errors and their dependence upon the patch spacing~\(\Delta\).

More quantitatively, in all cases, we found that the errors~\(\epsilon^{1,0}\) and~\(\epsilon^{1,1}\) followed the power law fit \({\epsilon}^{k_x,k_y} \approx 0.33 \cdot (0.7 \cdot \Delta)^p \).
The exponent is as expected for the various orders~\(p\) of inter-patch interpolation.
\cref{fig:iD2_iV2_eVsD_10_pN_m1dNE_Square_cD0p001_cV0p01_n10} is a `wosrt-case' because roundoff error started affecting the results for very small micro-grid spacing~\(\delta \lesssim 10^{-5}\) for the smallest patch ratio \(r=0.0001\) (circles).
The errors for wavenumber \((k_x,k_y)=(2,1)\) were found to behave similarly but with a larger coefficient, namely they follow the power law \({\epsilon}^{2,1} \approx 0.33 \cdot (1.33 \cdot \Delta)^p \).
Consequently, because the errors in the macroscale eigenvalues follow these power laws, \(\epsilon\propto\Delta^p\), \emph{the polynomial patch scheme has errors well controlled by both the order~\(p\) of inter-patch interpolation and the inter-patch spacing~\(\Delta\)}.

\section{Large computational savings in time simulations}
% ++++++++++++++++++++++++++++++++++++++++++++++++++++++++++++++++++++
\label{sec:cmptnlSvngs_gLnrWve}

This section quantifies and demonstrates the potentially large computational savings of the staggered patch schemes for wave systems.

Consider the staggered patch grid of \cref{fig:PtchGrd_n2t0} with patch spacing~\(\Delta\) and patch size~\(l\).
For patch ratio \(r := l/(2\Delta)\), it is straightforward to see that {the 2D staggered patch schemes compute only within a small fraction~\(3 r^2\) of the area of the full domain}.
For example, for \(r = 0.1, 0.01, 0.001, 0.0001\), the staggered patch schemes compute over the small fractions of area~\(0.03,\, 3 \cdot 10^{-4},\, 3 \cdot 10^{-6},\, 3 \cdot 10^{-8}\) respectively.
Similarly, in \(d\)~spatial dimensions, \(d\)-D staggered patch schemes would compute only within the small fraction~\((d+1) r^d\) of the volume~\(L^d\) of the full \(d\)-D domain.
Thus, a staggered patch scheme computes only within a small fraction of the space in the full domain and so we expect large computational savings.
But in a patch scheme, there is both the overhead of inter-patch interpolation and the more complicated pattern of memory access, and so we studied an example implementation to confirm \text{the potential savings.}

\begin{SCfigure}
    \caption{\label{fig:compTime_RU_PEI_gLnrWve_TpByTmd_vs_n}%
      The plotted point data are the ratio~\(T_p / T_\mu\) of the measured compute times of the staggered patch schemes (with {patch coupling via sparse matrix multiplication}) to that of the full-domain microscale model, for different orders of interpolation~\(p\), different patch ratios~\(r\), and across different numbers~\(n\) of sub-patch micro-grid intervals.
      Solid lines represent formula~\cref{eqn:TpByTm_gLnrWve} for~\(T_p / T_\mu\) using the estimated~\(T_M = 0.062\,\mu\)s and the respective estimated inter-patch coupling compute times~\(T_C\) for each patch scheme.
    }
    \includegraphics{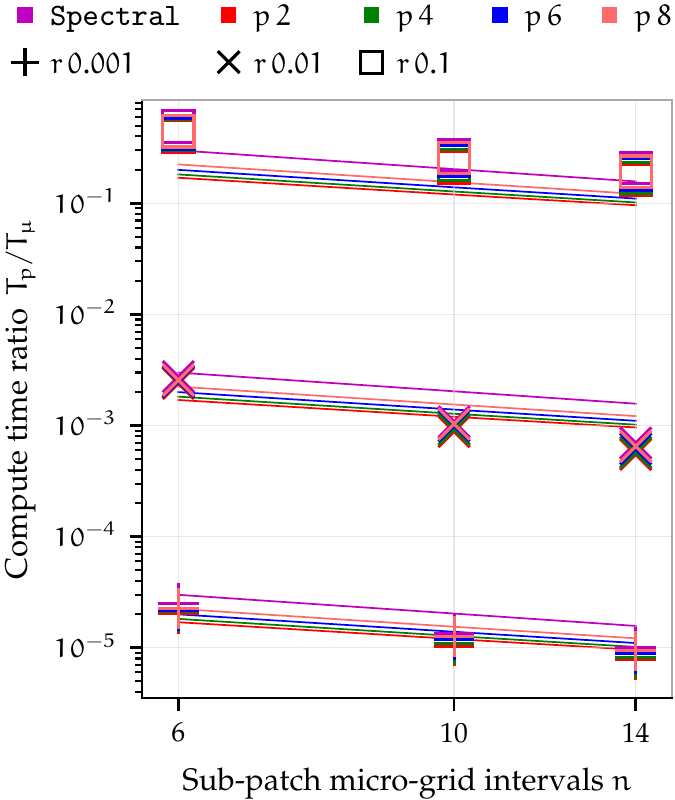}
\end{SCfigure}

\textcite[\S3.7.3]{Divahar2022_AcrteMltiscleSmltnOfWveLkeSystms}  details the derivation of the following formula for the crucial ratio of compute times for one evaluation of the time derivative~\eqref{eqn:dynSys_pN_gnrcWve}:
\begin{align}\label{eqn:TpByTm_gLnrWve}
    T_p/T_\mu = (T_C/T_M)r^2 \left( \frac{24}{n} - \frac{64}{3 n^2} \right)
      +  3 r^2 \left( 1 - \frac{16}{9 n} + \frac{8}{9 n^2}\right)
\end{align}
in which \(T_M\)~is the average compute time of the time derivative of one state variable in the microscale model~\cref{eqs:FDEs_mN_gLnrWve}; 
\(T_C\)~is the average compute time of one patch-edge value by the inter-patch coupling; \(r\)~is the patch ratio and \(n\)~is the number of sub-patch micro-grid intervals.
The crucial ratio~\(T_C/T_M\) in~\cref{eqn:TpByTm_gLnrWve} of coupling compute time to model compute time encapsulates both 
the details of the specific patch scheme (such as spectral or polynomial patch scheme, and interpolation order~\(p\));
and the details of the specific implementation (such as the algorithmic choices, data structures, serial, vector, or parallel computations).

\textcite[\S3.7.4]{Divahar2022_AcrteMltiscleSmltnOfWveLkeSystms} details the measurements of the compute time~\(T_\mu\) of the full domain model~\eqref{eqs:FDEs_mN_gLnrWve}, and measurements of the compute time~\(T_p\) of its various patch scheme implementations.
The compute times were measured on a workstation with Intel i7-6900k processor and \(64\)\textsc{gb} \textsc{ddr4~ram} using a specific implementation in Julia \parencite{Bezanson2017_JuliaAFrshAprchToNmrclCmptng}.
From the measurements of compute times~\(T_\mu\) and \(T_p\), the parameters~\(T_C\) and~\(T_M\) (for each patch scheme) of the expression~\cref{eqn:TpByTm_gLnrWve} were estimated as detailed by \textcite[\S3.7.4, Figs.~3.7.1 and 3.7.2]{Divahar2022_AcrteMltiscleSmltnOfWveLkeSystms}.
The ratio~\(T_p/T_\mu\) of the measured compute times in \cref{fig:compTime_RU_PEI_gLnrWve_TpByTmd_vs_n} shows that the reduction in compute time achieved by the patch scheme is proportional to the expected~\(r^2\), but with a coefficient somewhat larger than~\(3\) (and depending upon~\(n\)).

The main feature of the compute time ratio~\eqref{eqn:TpByTm_gLnrWve} is that both terms are proportional to~\(r^2\), the square of the patch ratio, as clearly evident in the computational measurements of \cref{fig:compTime_RU_PEI_gLnrWve_TpByTmd_vs_n}.
The coefficient of this~\(r^2\) behaviour is approximately \(3+\tfrac{24}n(T_C/T_M)\) as evident in the decrease of~\(T_p/T_\mu\) with~\(n\) in \cref{fig:compTime_RU_PEI_gLnrWve_TpByTmd_vs_n}.
This evidence demonstrates that the patch scheme may make accurate macroscale system-level predictions with speed-ups of up to a factor of~\(100\,000\) or more.

In problems with more involved microscale physics, \(T_M\)~will be larger, the relative compute cost of the inter-patch coupling will be smaller, and the corresponding potential speed-up could be even more.

\section{Conclusion}
% ++++++++++++++++++++++++++++++++++++++++++++++++++++++++++++++++++++

For large-scale waves (e.g., the planetary atmospheres, oceans, floods, and tsunamis), resolving the large range of spatial scales requires very many variables, leading to prohibitively high computational costs.
The equation-free multiscale modelling is a well-developed, powerful, and flexible approach to reducing the computational cost for dissipative systems.
But the small dissipation in waves poses a significant challenge to further developing the equation-free multiscale modelling methods, especially in multiple dimensions.
This article develops two novel families of equation-free multiscale 2D schemes, namely spectral (\cref{ssc:spctrlCplng}) and polynomial (\cref{ssc:plynmlCplng}) staggered patch schemes.
Qualitative exploration of the patch scheme eigenvalues (\cref{sec:ptchSchmsAreAcrte_gLnrWve}) shows both the structure and accuracy of the eigenvalues of macroscale modes.

A study of sensitivity to roundoff errors (\cref{sec:ptchSchmsAreNtSnstveToNmrclErrs_gLnrWve}) establishes the robustness of the developed staggered patch schemes.
This insensitivity to roundoff errors empowers using the multiscale schemes over a wide range of length scales and allows using advanced optimisations in finite precision computing.

Comprehensive eigenvalue analysis (\cref{sec:ptchSchmsAreStble_gLnrWve}) over a wide range of parameters shows that the developed patch schemes are stable.
Characterising the dependence of eigenvalues errors (\cref{{sec:ptchSchmsAreCnsstnt_gLnrWve}}) on the patch spacing~\(\Delta\) establishes that staggered patch schemes are consistent with the given microscale model.
Specifically, the spectral patch scheme is uniformly consistent without any dependence on the patch spacing, whereas the polynomial patch schemes are consistent to the same order~\(p\) of interpolation with decreasing patch spacing~\(\Delta\) (i.e., macroscale errors decrease as \(\Delta^p\)).

Theoretical quantification of the computational complexity and the measured compute times of the multiscale staggered patch schemes agree (\cref{sec:cmptnlSvngs_gLnrWve}).
Both demonstrate that the patch schemes may make accurate macroscale predictions with speed-ups of up to a factor of~\(100\,000\) or more for the weakly damped linear waves.
Compared to the considered simple weakly damped linear waves, modelling more complex physical processes leads to large model compute time~\(T_M\), compared to the inter-patch coupling compute time~\(T_C\), hence smaller~\(T_C/T_M\) in \cref{eqn:TpByTm_gLnrWve}, and so smaller~\(T_p/T_\mu\), leading to much larger potential speed-ups.
All the demonstrated computational speed-ups in this article are for a 2D spatial domain: larger speed-ups are feasible in more spatial dimensions.
Thus this work provides the essential foundation for efficient large-scale simulation of challenging nonlinear multiscale waves.

\paragraph*{Acknowledgments}
Parts of this research were supported by the Australian Research Council grants DP150102385 and DP200103097.
The work of \textsc{I.G.K.} was partially supported by a \textsc{muri} grant by the US Army Research Office (Drs.~S. Stanton and M. Munson).
J. Divahar was supported by an Australian Government Research Training Program (RTP) Scholarship.

\end{document}